\documentclass[12pt]{article}
\usepackage{amsfonts}
\usepackage{amsmath}
\usepackage{hyperref}

\def\squarebox#1{\hbox to #1{\hfill\vbox to #1{\vfill}}}
\newcommand{\qed}{\hspace*{\fill}
\vbox{\hrule\hbox{\vrule\squarebox{.667em}\vrule}\hrule}\smallskip}
\newtheorem{teorema}{Theorem}[section]
\newtheorem{lema}[teorema]{Lemma}
\newtheorem{corolario}[teorema]{Corollary}
\newtheorem{proposicao}[teorema]{Proposition}
\newtheorem{defi}[teorema]{Definition}
\newenvironment{profe}{\noindent {\bf Proof:}}{\hfill $\qed $ \newline}
\begin{document}

\title{The Isotropy Representation of a Real Flag Manifold: Split Real Forms}
\author{Mauro Patr\~{a}o \thanks{%
Department of Mathematics. Universidade Federal de Brasilia. Bras\'{\i}lia -
DF, Brazil. Supported by CNPq grant no.\ 310790/09-3} \and Luiz A. B. San
Martin \thanks{%
Institute of Mathematics. Universidade Estadual de Campinas. Campinas - SP,
Brazil. Supported by CNPq grant no.\ 303755/09-1, FAPESP grant no.\
2012/18780-0 and CNPq/Universal grant no 476024/2012-9.}}
\date{}
\maketitle

\begin{abstract}
We study the isotropy representation of real flag manifolds associated to
simple Lie algebras that are split real forms of complex simple Lie
algebras. For each Dynkin diagram the invariant irreducible subspaces for
the compact part of the isotropy subgroup are described. Contrary to the
complex flag manifolds the decomposition into irreducible components is not
in general unique, since there are cases with infinitely many invariant
subspaces.
\end{abstract}

\noindent \textit{AMS 2010 subject classification}: Primary: 14M15,
Secondary: 57R22.

\noindent \textit{Key words an phrases}. Flag manifold (generalized),
Isotropy representation,  Real semi-simple Lie groups, Split-real forms of
complex semi-simple Lie algebras.

\section{Introduction\label{secintro}}

This paper studies the isotropy representation of (generalized) real flag
manifolds associated to a noncompact real simple Lie algebra $\mathfrak{g}$.
Here we consider the case where $\mathfrak{g}$ is a split real form of a
complex simple Lie algebra.

A flag manifold of $\mathfrak{g}$ is a coset space $\mathbb{F}_{\Theta
}=G/P_{\Theta }$ where $G$ is any connected Lie group with Lie algebra $%
\mathfrak{g}$ and $P_{\Theta }\subset G$ is a parabolic subgroup. The Lie
algebra $\mathfrak{p}_{\Theta }$ of $P_{\Theta }$ is a parabolic subalgebra
which is the sum of the eigenspaces of the nonnegative eigenvalues of $%
\mathrm{ad}\left( H_{\Theta }\right) $ with $H_{\Theta }\in \mathfrak{g}$ a
suitable chosen element. If $K\subset G$ is a maximal compact subgroup and $%
K_{\Theta }=K\cap P_{\Theta }$ then $\mathbb{F}_{\Theta }=K/K_{\Theta }$ as
well.

The two presentations $\mathbb{F}_{\Theta }=G/P_{\Theta }$ and $\mathbb{F}%
_{\Theta }=K/K_{\Theta }$ yield the isotropy representations of $P_{\Theta }$
and $K_{\Theta }$ on the tangent space $T_{b_{\Theta }}\mathbb{F}_{\Theta }$
at the origin $b_{\Theta }$. The $K_{\Theta }$-representation is obtained by
restricting the $P_{\Theta }$-representation.

Our objective in this paper is to describe the isotropy representation of $%
K_{\Theta }$. This means that the invariant and irreducible subspaces of $%
T_{b_{\Theta }}\mathbb{F}_{\Theta }$ must be obtained as well as the
possible decompositions 
\begin{equation}
T_{b_{\Theta }}\mathbb{F}_{\Theta }=V_{1}\oplus \cdots \oplus V_{k}
\label{fordecintro}
\end{equation}%
into $K_{\Theta }$-invariant irreducible components.

The description of the isotropy representatin of $K_{\Theta }$ is essential
to get $K$-invariant geometries on $\mathbb{F}_{\Theta }$. For example the $%
K $-invariant Riemannian metrics on $\mathbb{F}$ are given by the $K_{\Theta
}$-invariant inner products on $T_{b_{\Theta }}\mathbb{F}_{\Theta }$, which
in turn are direct sum of invariant inner products on the components of a
decomposition (\ref{fordecintro}).

Too look at the $K_{\Theta }$-representation we consider first the the
isotropy representation of $P_{\Theta }$. It is completely determined by the
restriction to its Levi component $Z_{\Theta }$, which is the centralizer in 
$G$ of $H_{\Theta }$. The group $Z_{\Theta }$ is reductive, so that its
representation decomposes as a sum of invariant irreducible subspaces. This
decomposition is unique and coincide with the decomposition for the ensuing
representation of the Lie algebra $\mathfrak{z}_{\Theta }$ of $Z_{\Theta }$.
In fact, each $\mathfrak{z}_{\Theta }$-irreducible component is a sum of
root spaces (for a Cartan subalgebra) associated to different roots for
different componets. This implies uniqueness of the decomposition. For the
moment we write the $\mathfrak{z}_{\Theta }$-decomposition as 
\begin{equation*}
T_{b_{\Theta }}\mathbb{F}_{\Theta }=W_{1}^{\mathfrak{z}}\oplus \cdots \oplus
W_{n}^{\mathfrak{z}}.
\end{equation*}

The subspaces $W_{i}^{\mathfrak{z}}$ are invariant by $K_{\Theta }$ since $%
K_{\Theta }\subset Z_{\Theta }$. Hence we are faced to the following
problems:

\begin{enumerate}
\item Find the $K_{\Theta }$-invariant irreducible subspaces inside each $%
W_{i}^{\mathfrak{z}}$. This includes the question of deciding whether $%
W_{i}^{\mathfrak{z}}$ is $K_{\Theta }$-irreducible.

\item Among the invariant subspaces of item (1), find those pairs $U_{1}$, $%
U_{2}$ such that the $K_{\Theta }$-representations on them are equivalent.
Given such a pair we get further invariant irreducible subspaces contained
in $U_{1}\oplus U_{2}$ as graphs of operators $T:U_{1}\rightarrow U_{2}$,
intertwining the representations on $U_{1}$ and $U_{2}$.
\end{enumerate}

The answers to these two questions give the full picture of the $K_{\Theta }$%
-invariant subspaces.

At this point it is worthwhile to compare the real flag manifolds with the
complex ones. In the complex case the above questions have trivial answers:
The subspaces $W_{1}^{\mathfrak{z}}$ are $K_{\Theta }$-irreducible and no
two of them are equivalent. This is due to the fact that in a complex Lie
group $K_{\Theta }$ is a compact real form of the semi-simple component $%
G\left( \Theta \right) $ of $Z_{\Theta }$, which is also a complex group. So
that the equivalence classes of $K_{\Theta }$-representations are in
bijection to the $G\left( \Theta \right) $-representations.

On the contrary for real flag manifolds new phenomena occur: There are $%
\mathfrak{z}_{\Theta }$-irreducible subspaces that are not $K_{\Theta }$%
-irreducible and there are equivalent $K_{\Theta }$-invariant irreducible
subspaces. Such equivalence gives rise to continuous sets of invariant
subspaces and to the nonuniqueness of the decompositions (\ref{fordecintro}).

The basic differences of the real case to the complex one is that $K$ is not
in general simple and $K_{\Theta }$ is not connected (if $\mathfrak{g}$ is a
split real form). When $K$ is not simple we get a supply of $K_{\Theta }$%
-invariant subspaces as tangent spaces to the orbits through the origin $%
b_{\Theta }$ of the simple components of $K$. In many cases these tangent
spaces decompose the $\mathfrak{z}_{\Theta }$-irreducible subspaces. The
fact that $K_{\Theta }$ is not connected requires a separate analysis of the
representations of its group of connected components, the so called $M$%
-group.

Now we describe the contents of the paper. Section \ref{secisotro} contains
generalities about isotropy representations.

The main technical part of the paper starts at Section \ref{secmequiv} where
we look at the representations of the discrete group $M$. This is the
centralizer in $K$ of the Cartan subalgebra $\mathfrak{a}$ and contains
information about the group of connected components of any $K_{\Theta }$.
Also $M=K_{\Theta }$ if $\mathbb{F}_{\Theta }$ is the maximal flag manifold.
The one dimensional root spaces $\mathfrak{g}_{\alpha }$ are $M$-invariant
thus defining representations of $M$. For the roots $\alpha $ and $\beta $
we put $\alpha \sim _{M}\beta $ if the representations of $M$ on $\mathfrak{g%
}_{\alpha }$ and $\mathfrak{g}_{\beta }$ are equivalent. The purpose of
Section \ref{secmequiv} is to find $M$-equivalence classes of roots. After
some preparations we proceed to a case by case analysis of the diagrams. For
each case the $M$-equivalence classes are described at the beginning of the
corresponding subsection. For the classical diagrams there are exceptions
since in low dimension the sizes of the classes tend to increase. The
detemination of the $M$-equivalence classes furnishes the complete picture
of the isotropy representation on the maximal flag manifolds. They will be
also a basic tool to detect inequivalent subrepresentations in the other
flag manifolds.

Section \ref{seclemaux} is preparatory. There we prove several lemmas to be
applied in the study of isotropy representations on the partial flag
manifolds. Some of these lemmas have independent interest, like Lemma \ref%
{lemtransimplelong} which ensures transitivity of the Weyl group on the set
of weights of a given representation. This fact is far from to be true for
general representations.

In Section \ref{secirreduc} we go through the isotropy representations of
the partial flag manifolds, again in a case by case analysis. For the
classical diagrams we use their standard realizations as algebras of
matrices: $A_{l}=\mathfrak{sl}\left( l+1,\mathbb{R}\right) $, $B_{l}=%
\mathfrak{so}\left( l,l+1\right) $, $C_{l}=\mathfrak{sp}\left( l,\mathbb{R}%
\right) $ and $D_{l}=\mathfrak{so}\left( l,l\right) $. These realizations
allow the use of nice expressions for the roots. The analysis of the
classical diagrams have the following pattern: First we describe the $%
\mathfrak{z}_{\Theta }$-irreducible components. Then we check their $%
K_{\Theta }$-irreducibility and finally we look at equivalence between
irreducible subspaces. The results are summarized at the end of each
corresponding subsection. Regarding to the exceptional diagrams, $G_{2}$ is
clear by its low dimensionality. For $E_{6}$, $E_{7}$ and $E_{8}$, it
follows easily by the general lemmas of Section \ref{seclemaux} that the $%
K_{\Theta }$-invariant subspaces are the $\mathfrak{z}_{\Theta }$%
-irreducible components. \ As to $F_{4}$ we refrain to make a detailed and
annoying description of the fifteen flag manifolds. Besides the maximal flag
manifold, where the picture is given by the $M$-equivalence classes, we just
look at a minimal flag manifold.

In conclusion we say that our initial motivation to study the isotropy
representation came from the attempt to understand the $K$-invariant
Riemannian metrics on the real flag manifolds. There is an extensive
literature on invariant Riemannian geometry on complex flag manifolds. See
for example Burstall-Rawnsley \cite{br}, Burstall-Salamon \cite{bs},
Negreiros \cite{neg}, San Martin-Negreiros \cite{smneg}, San Martin-Silva 
\cite{smrit}, and Wang-Ziller \cite{wz}, and references therein. In a
complex flag manifold the isotropy representation has a unique decomposition
into invariant irreducible components, which makes the set of invariant
Riemannian metrics a finite dimensional manifold. Our results in this paper
show the existence of infinitely many decompositions on a real flag
manifold, pointing to a great richness of the invariant Riemannian geometry.

\section{Isotropy representation\label{secisotro}}

Let $\mathfrak{g}$ be a split real form of a complex simple Lie algebra, $%
\mathfrak{g}=\mathfrak{k}\oplus \mathfrak{s}$ be a Cartan decomposition and $%
\mathfrak{a}\subset \mathfrak{s}$ be a maximal abelian subalgebra. Denote by 
$\Pi $ the associated set of roots and by 
\begin{equation*}
\mathfrak{g}=\mathfrak{g}_{0}\oplus \sum_{\alpha \in \Pi }\mathfrak{g}%
_{\alpha }
\end{equation*}%
the associated root space decomposition. Denote by $G$ the group of inner
automorphisms of $\mathfrak{g}$, which is the subgroup of $\mathrm{Gl}(%
\mathfrak{g})$ generated by $\exp \mathrm{ad}(\mathfrak{g})$. Let $K$ be the
subgroup of $G$ generated by $\mathrm{ad}(\mathfrak{k})$. Fixing a set $\Pi
^{+}$ of positive roots let $\Sigma $ be the corresponding set of simple
roots. We denote by $\mathfrak{a}^{+}=\{H\in \mathfrak{a}:\forall \alpha \in
\Sigma $, $\alpha \left( H\right) >0\}$ the Weyl chamber associated to $%
\Sigma $.

A subset $\Theta \subset \Sigma $ defines the parabolic subalgebra of type $%
\Theta $ given by 
\begin{equation*}
\mathfrak{p}_{\Theta }=\mathfrak{g}_{0}\oplus \sum_{\alpha \in \Pi ^{+}}%
\mathfrak{g}_{\alpha }\oplus \sum_{\alpha \in \langle \Theta \rangle ^{-}}%
\mathfrak{g}_{\alpha },
\end{equation*}%
where $\langle \Theta \rangle ^{-}$ is the set of negative roots generated
by $\Theta $. The standard parabolic subgroup $P_{\Theta }$ defined by $%
\Theta $ is the normalizer of $\mathfrak{p}_{\Theta }$ in $G$. The
associated flag manifold is defined by $\mathbb{F}_{\Theta }=G/P_{\Theta }$.
Since $K$ acts transitively on $\mathbb{F}_{\Theta }$, this flag manifold
can be given by $\mathbb{F}_{\Theta }=K/K_{\Theta }$, where $K_{\Theta
}=P_{\Theta }\cap K$.

When $\Theta =\emptyset $ we get the minimal parabolic subalgebra $\mathfrak{%
p}=\mathfrak{p}_{\emptyset }$. In this case the subscript is omited and the
maximal flag manifold is written $\mathbb{F}=G/P$. We have $\mathbb{F}=K/M$,
where $M=K_{\emptyset }$ is the centralizer of $\mathfrak{a}$ in $K$.

For an alternative description of the parabolic subalgebra write 
\begin{equation*}
\mathfrak{a}_{\Theta }=\{H\in \mathfrak{a}:\forall \alpha \in \Theta ,\alpha
\left( H\right) =0\}
\end{equation*}%
for the anihilator of $\Theta $. Let $H_{\Theta }$ be characteristic for $%
\Theta $, that is $H_{\Theta }$ is in the \textquotedblleft partial
chamber\textquotedblright\ $\mathfrak{a}_{\Theta }\cap \mathrm{cl}\mathfrak{a%
}^{+}$ and satisfies 
\begin{equation*}
\Theta =\{\alpha \in \Sigma :\alpha (H_{\Theta })=0\}.
\end{equation*}%
Then 
\begin{equation*}
\mathfrak{p}_{\Theta }=\sum_{\lambda \geq 0}V_{\lambda }\left( H_{\Theta
}\right)
\end{equation*}%
where $V_{\lambda }\left( H_{\Theta }\right) =\sum_{\alpha \left( H_{\Theta
}\right) =\lambda }\mathfrak{g}_{\alpha }$ is the $\lambda $-eigenspace of $%
\mathrm{ad}\left( H_{\Theta }\right) $. Clearly any $H_{\Theta }$ satisfying
(\ref{forcaracteris}) yield the same $\mathfrak{p}_{\Theta }$, although the
eigenspaces $V_{\lambda }\left( H_{\Theta }\right) $ may change.

The centralizer of $H_{\Theta }$, $\mathfrak{z}_{\Theta }=\mathrm{cent}_{%
\mathfrak{g}}\left( H_{\Theta }\right) =\sum_{\alpha \left( H_{\Theta
}\right) =0}\mathfrak{g}_{\alpha }$ is the Levi component of $\mathfrak{p}%
_{\Theta }$. It is a reductive Lie algebra that decomposes as 
\begin{equation*}
\mathfrak{z}_{\Theta }=\mathfrak{g}\left( \Theta \right) \oplus \mathfrak{a}%
_{\Theta }
\end{equation*}%
where the semi-simple component $\mathfrak{g}\left( \Theta \right) $ is the
subalgebra generated by $\mathfrak{g}_{\alpha }$, $\alpha \in \pm \Theta $.
Since $\mathfrak{g}$ is a split real form, it follows that $\mathfrak{g}%
\left( \Theta \right) $ is also a split real form, having Cartan subalgebra
the subspace $\mathfrak{a}\left( \Theta \right) $ spanned by $H_{\alpha }$, $%
\alpha \in \Theta $ (where $\alpha \left( \cdot \right) =\langle H_{\alpha
},\cdot \rangle $). Put $G\left( \Theta \right) =\langle \exp \mathfrak{g}%
\left( \Theta \right) \rangle $ for the connected subgroup with Lie algebra $%
\mathfrak{g}\left( \Theta \right) $.

With this notation we have that $K_{\Theta }=\mathrm{Cent}_{K}\left(
H_{\Theta }\right) $ is the centralizer of $H_{\Theta }$ in $K$ and its Lie
algebra $\mathfrak{k}_{\Theta }=\mathrm{Cent}_{\mathfrak{k}}\left( H_{\Theta
}\right) =\mathfrak{z}_{\Theta }\cap \mathfrak{k}$. Also, $\left( K_{\Theta
}\right) _{0}\subset G\left( \Theta \right) $.

The nilpotent subalgebra 
\begin{equation*}
\mathfrak{n}_{\Theta }^{-}=\sum_{\alpha \in \Pi ^{-}\backslash \langle
\Theta \rangle ^{-}}\mathfrak{g}_{\alpha }=\sum_{\lambda <0}V_{\lambda
}\left( H_{\Theta }\right)
\end{equation*}%
complements $\mathfrak{p}_{\Theta }$ in $\mathfrak{g}$. Hence we identify $%
\mathfrak{n}_{\Theta }^{-}$ with the tangent space $T_{b_{\Theta }}\mathbb{F}%
_{\Theta }$ at the origin $b_{\Theta }$. Under this identification the
isotropy representations of $K_{\Theta }$ and $G\left( \Theta \right) $ are
just the adjoint representation, since $\mathfrak{n}_{\Theta }^{-}$ is
normalized by these groups. The same statement holds for the representations
of the Lie algebras $\mathfrak{k}_{\Theta }$, $\mathfrak{g}\left( \Theta
\right) $ and $\mathfrak{z}_{\Theta }$.

Since $\mathfrak{z}_{\Theta }$ is reductive its representation on $\mathfrak{%
n}_{\Theta }$ is a direct sum 
\begin{equation*}
\mathfrak{n}_{\Theta }=\sum_{\sigma }V_{\Theta }^{\sigma }
\end{equation*}%
where the subspaces $V_{\Theta }^{\sigma }$ are $\mathfrak{z}_{\Theta }$%
-invariant and irreducible. Here we use $\sigma $ to distinguish the
different invariant subspaces.

\begin{proposicao}
\label{propcomproot}Each $\mathfrak{z}_{\Theta }$-invariant and irreducible
subspace $V_{\Theta }^{\sigma }$ is a direct sum of root spaces, 
\begin{equation*}
V_{\Theta }^{\sigma }=\sum \mathfrak{g}_{\alpha }
\end{equation*}%
where the sum extended to a subset of roots $\Pi _{\Theta }^{\sigma }\subset
\Pi ^{-}\backslash \langle \Theta \rangle ^{-}$. Conversely if $\alpha \in
\Pi ^{-}\backslash \langle \Theta \rangle ^{-}$ then $\mathfrak{g}_{\alpha }$
is contained in a unique $\mathfrak{z}_{\Theta }$-component denoted by $%
V_{\Theta }\left( \alpha \right) $. We write $\Pi _{\Theta }\left( \alpha
\right) $ for the roots $\beta $ with $\mathfrak{g}_{\beta }\subset
V_{\Theta }\left( \alpha \right) $.
\end{proposicao}

\begin{profe}
This follows by a standard argument using the fact that $\mathfrak{a}\subset 
\mathfrak{z}_{\Theta }$. In fact, if $V$ is a $\mathfrak{z}_{\Theta }$%
-invariant subspace and $X=\sum X_{\alpha }\in V$ then 
\begin{equation*}
\mathrm{ad}\left( H\right) X=\sum \alpha \left( H\right) X_{\alpha }\in V
\end{equation*}%
if $H\in \mathfrak{a}$. By taking suitable values of $H\in \mathfrak{a}$ one
concludes that each component $X_{\alpha }\in V$, so that $\mathfrak{g}%
_{\alpha }\subset V$. The last statement follows directly from the fact tha $%
\mathfrak{n}_{\Theta }^{-}$ is the direct sum of the roots spaces as well as
the $\mathfrak{z}_{\Theta }$-components.
\end{profe}

The restriction to $\mathfrak{g}\left( \Theta \right) $ of the $\mathfrak{z}%
_{\Theta }$-representation on $V_{\Theta }^{\sigma }$ is also irreductible.
This is because $\mathfrak{z}_{\Theta }=\mathfrak{g}\left( \Theta \right)
\oplus \mathfrak{a}_{\Theta }$ with $\mathfrak{a}_{\Theta }$ the center of $%
\mathfrak{g}\left( \Theta \right) $, so that $\mathrm{ad}\left( H\right) $
is a scalar $\lambda \cdot \mathrm{id}$ in $V_{\Theta }^{\sigma }$ for any $%
H\in \mathfrak{a}_{\Theta }$. Hence, any $\mathfrak{g}\left( \Theta \right) $%
-invariant subspace $U\subset V_{\Theta }^{\sigma }$ is also $\mathfrak{z}%
_{\Theta }$-invariant, ensuring that $V_{\Theta }^{\sigma }$ is $\mathfrak{g}%
\left( \Theta \right) $-irreducible.

The weight spaces of the representation of $\mathfrak{g}\left( \Theta
\right) $, w.r.t. $\mathfrak{a}\left( \Theta \right) $, are root spaces of $%
\mathfrak{g}$, so that the weights of the representation are restrictions to 
$\mathfrak{a}\left( \Theta \right) $ of some roots $\alpha \in \Pi
^{-}\backslash \langle \Theta \rangle ^{-}$. There is just one highest
weight, say $\mu _{\sigma }$, and two representations of $\mathfrak{g}\left(
\Theta \right) $ on $V_{\Theta }^{\sigma _{1}}$ and $V_{\Theta }^{\sigma
_{2}}$ are equivalent if and only if $\mu _{\sigma _{1}}=\mu _{\sigma _{2}}$%
. (We note that different representations of $\mathfrak{z}_{\Theta }$ on $%
V_{\Theta }^{\sigma _{1}}$ and $V_{\Theta }^{\sigma _{2}}$ cannot be
equivalent, even if the $\mathfrak{g}\left( \Theta \right) $-representations
are equivalent.)

The subspaces $V_{\Theta }^{\sigma }$ are also invariant and irreductible by 
$G\left( \Theta \right) $, since by definition this group is connected.
Hence $V_{\Theta }^{\sigma }$ is invariant by the identity component $\left(
K_{\Theta }\right) _{0}$ of $K_{\Theta }$, because $\left( K_{\Theta
}\right) _{0}$ $\subset G\left( \Theta \right) $.\ As to $K_{\Theta }$ we
have $K_{\Theta }=M\cdot \left( K_{\Theta }\right) _{0}$ which ensures that $%
V_{\Theta }^{\sigma }$ is $K_{\Theta }$-invariant, because $M$ leaves
invariant each root space.

Our objective is to get the invariant irreducible subspaces of $\mathfrak{n}%
_{\Theta }^{-}$ by the $K_{\Theta }$ representation, which is equivalent to
the isotropy representation of the flag $\mathbb{F}_{\Theta }$.

In view of the above discussion we are reduced to the following questions:

\begin{enumerate}
\item Describe the irreducible components $V_{\Theta }^{\sigma }$ of the $%
\mathfrak{z}_{\Theta }$ representation.

\item Find the $K_{\Theta }$-invariant subspaces of each $V_{\Theta
}^{\sigma }$.

\item Find the pairs of irreducible subspaces having equivalent $K_{\Theta }$%
-representations.
\end{enumerate}

Finally we note that if $H_{\Theta }\in \mathfrak{a}_{\Theta }\cap \mathrm{cl%
}\mathfrak{a}^{+}$ then an eigenspace $V_{\lambda }\left( H_{\Theta }\right)
=\sum_{\alpha \left( H_{\Theta }\right) =\lambda }\mathfrak{g}_{\alpha }$ is
contained in $\mathfrak{n}_{\Theta }^{-}$ if $\lambda <0$ and is invariant
by $\mathfrak{z}_{\Theta }$. Hence $V_{\lambda }\left( H_{\Theta }\right) $
is the direct sum of some irreducible components $V_{\Theta }^{\sigma }$.
This remark will be used later to determine the irreducible components $%
V_{\Theta }^{\sigma }$. Actually, in some cases an eigenspace $V_{\lambda
}\left( H_{\Theta }\right) $ is irreducible and hence is itself a component.

\section{$M$-equivalence classes\label{secmequiv}}

Let \ $M=\mathrm{Cent}_{K}\left( \mathfrak{a}\right) $ be the centralizer of 
$\mathfrak{a}$ in $K$. It is known that $M\subset K_{\Theta }=M(K_{\Theta
})_{0}$. Also, any root space $\mathfrak{g}_{\alpha }$ is $M$-invariant. In
this section we determine the pairs of root spaces $\mathfrak{g}_{\alpha }$, 
$\mathfrak{g}_{\beta }$ having equivalent representations of $M$. This will
be used later to check equivalence or nonequivalence of $K_{\Theta }$%
-representations on invariant subspaces.

\begin{defi}
The roots $\alpha $ and $\beta $ are said to be $M$-equivalent (in symbols $%
\alpha \sim _{M}\beta $) if the representations of $M$ on $\mathfrak{g}%
_{\alpha }$ and $\mathfrak{g}_{\beta }$ are equivalent. We write $\left[
\alpha \right] _{M}$ for the $M$-equivalence class of the root $\alpha $.
\end{defi}

If $\mathfrak{g}$ is a split real form of a complex semi-simple Lie algebra
then $M$ is a discrete abelian subgroup equals to 
\begin{equation*}
M=\{m_{\gamma }=\exp (\pi iH_{\gamma }^{\vee }):\gamma \in \Pi \}
\end{equation*}%
where $H_{\gamma }^{\vee }=\frac{2H_{\gamma }}{\langle \gamma ,\gamma
\rangle }$ is the co-root associated to $\gamma $ and $H_{\gamma }$ is
defined by $\gamma (H)=\langle H_{\gamma },H\rangle $, $H\in \mathfrak{a}$.
In the above formula the exponential $\exp (\pi iH_{\gamma }^{\vee })$ is in
the complex group $\mathrm{Aut}\left( \mathfrak{g}_{\mathbb{C}}\right) $,
where $\mathfrak{g}_{\mathbb{C}}$ is the complexification of $\mathfrak{g}$
(see \cite{knp}, Theorems 7.53 and 7.55).

The following statement gives a necessary and sufficient condition for the $%
M $-equivalence between the roots $\alpha $ and $\beta $.

\begin{proposicao}
\label{propmequivalent}The root $\alpha $ and $\beta $ are $M$-equivalent if
and only if, for every $\gamma \in \Pi $ we have 
\begin{equation*}
\frac{2\langle \gamma ,\alpha \rangle }{\langle \gamma ,\gamma \rangle }%
\equiv \frac{2\langle \gamma ,\beta \rangle }{\langle \gamma ,\gamma \rangle 
} \quad \mbox{mod}2.
\end{equation*}
\end{proposicao}

\begin{profe}
Take a root $\gamma $ and write as above $m_{\gamma }=\exp (\pi iH_{\gamma
}^{\vee })$. If $X\in \mathfrak{g}_{\alpha }$ and $Y\in \mathfrak{g}_{\beta
} $ then 
\begin{equation*}
\mathrm{Ad}(m_{\gamma })X=e^{\pi i\alpha (H_{\gamma }^{\vee })}X\qquad %
\mbox{and}\qquad \mathrm{Ad}(m_{\gamma })Y=e^{\pi i\beta (H_{\gamma }^{\vee
})}Y,
\end{equation*}%
by definition of $m_{\gamma }$. It follows that $\alpha \sim _{M}\beta $ if
and only if $e^{\pi i\alpha (H_{\gamma }^{\vee })}=e^{\pi i\beta (H_{\gamma
}^{\vee })}$, which is equivalent to 
\begin{equation*}
\frac{2\langle \gamma ,\alpha \rangle }{\langle \gamma ,\gamma \rangle }%
\equiv \frac{2\langle \gamma ,\beta \rangle }{\langle \gamma ,\gamma \rangle 
} \quad \mbox{mod}2
\end{equation*}%
as desired.
\end{profe}

As a corollary we get the following necessary condition.

\begin{corolario}
\label{cormequivalent}If $\alpha \sim _{M}\beta $ then $\langle \alpha
,\beta \rangle =0$.
\end{corolario}

\begin{profe}
Suppose that $\langle \alpha ,\beta \rangle \neq 0$. Then we have the
following possibilities for the Killing numbers:

\begin{enumerate}
\item $\alpha $ and $\beta $ have the same length and the angle between them
is $60%
{{}^\circ}%
$ or $120%
{{}^\circ}%
$. Then 
\begin{equation*}
\frac{2\langle \alpha ,\alpha \rangle }{\langle \alpha ,\alpha \rangle }%
=2\qquad \mbox{and}\qquad \frac{2\langle \alpha ,\beta \rangle }{\langle
\alpha ,\alpha \rangle }=\pm 1
\end{equation*}%
showing that $\alpha $ and $\beta $ are not $M$-equivalent.

\item The angle between $\alpha $ and $\beta $ is $45%
{{}^\circ}%
$ or $135%
{{}^\circ}%
$. If $\alpha $ is the long root then 
\begin{equation*}
\frac{2\langle \alpha ,\alpha \rangle }{\langle \alpha ,\alpha \rangle }%
=2\qquad \mbox{and}\qquad \frac{2\langle \alpha ,\beta \rangle }{\langle
\alpha ,\alpha \rangle }=\pm 1
\end{equation*}%
and $\alpha $ and $\beta $ are not $M$-equivalent. If otherwise $\beta $ is
the long root then we interchange the roles of $\alpha $ and $\beta $ to get
the same result.

\item The angle between $\alpha $ and $\beta $ is $30%
{{}^\circ}%
$ or $150%
{{}^\circ}%
$. Then 
\begin{equation*}
\frac{2\langle \alpha ,\alpha \rangle }{\langle \alpha ,\alpha \rangle }%
=2\qquad \mbox{and}\qquad \frac{2\langle \alpha ,\beta \rangle }{\langle
\alpha ,\alpha \rangle }=\pm 1,\pm 3,
\end{equation*}%
concluding the proof.
\end{enumerate}
\end{profe}

In the sequel we apply the above proposition and its corollary to find for
each Dynkin diagram the classes of $M$-equivalence between roots. The
following simple remarks are used throughout with no further reference.

\begin{enumerate}
\item $\alpha \sim _{M}\left( -\alpha \right) $. In fact the Cartan
involution $\theta $ is an equivalence between the $M$-representations in $%
\mathfrak{g}_{\alpha }$ and $\mathfrak{g}_{-\alpha }$ because $\theta \left(
m\right) =m$ if $m\in M$. This implies that $\mathrm{Ad}\left( m\right)
\circ \theta =\theta \circ \mathrm{Ad}\left( m\right) $ if $m\in M$ and
since $\theta \left( \mathfrak{g}_{\alpha }\right) =\mathfrak{g}_{-\alpha }$
the equivalence follows. Hence we are reduced to check $M$-equivalence
between positive roots alone.

\item The criterion of Proposition \ref{propmequivalent} implies easily that
if $w\in \mathcal{W}$ then $\alpha \sim _{M}\beta $ if and only if $w\alpha
\sim _{M}w\beta $. Hence it will will be enough to get the $M$-equivalence
classes for just one element in each orbit of the Weyl group, that is, for
one root in the simply laced diagrams and for one long root and a short root
in diagrams with multiple edges.
\end{enumerate}

We proceed now to look at the $M$-equivalences for each Dynkin diagram.

\subsection{Diagram $A_{l}$, $l\geq 1$}

We use the standard realization of $A_{l}$ where the positive roots are
written as $\lambda _{i}-\lambda _{j}$, $1\leq i<j\leq l+1$. There are two
cases:

\subsubsection{$A_{l}$, $l\neq 3$}

The classes of \textbf{\ }$M$-equivalence on the positive roots are
singletons. (That is the $M$-representation on different root spaces are not
equivalent.)

Since the Weyl group is transitive on the set of roots it is enough to fix a
specific root $\alpha $ and check that any positive root $\beta \neq \alpha $
is not $M$-equivalent to $\alpha $.

Suppose that $l>3$ and take $\alpha =\lambda _{1}-\lambda _{2}$. The
positive roots orthogonal to $\alpha $ are $\lambda _{i}-\lambda _{j}$ with $%
3\leq i<j$. By Corollary \ref{cormequivalent} we are reduced to check that
these roots are not $M$-equivalent to $\alpha =\lambda _{1}-\lambda _{2}$.
There are the following cases for $3\leq i<j$:

\begin{enumerate}
\item If $j<l+1$ then $\langle \gamma ,\lambda _{i}-\lambda _{j}\rangle \neq
0$ but $\langle \gamma ,\lambda _{1}-\lambda _{2}\rangle =0$ where $\gamma
=\lambda _{j}-\lambda _{j+1}$. Hence 
\begin{equation*}
\frac{2\langle \gamma ,\lambda _{1}-\lambda _{2}\rangle }{\langle \gamma
,\gamma \rangle }=0\qquad \mathrm{and}\qquad \frac{2\langle \gamma ,\lambda
_{i}-\lambda _{j}\rangle }{\langle \gamma ,\gamma \rangle }=\pm 1
\end{equation*}%
so that $\lambda _{1}-\lambda _{2}$ is not $M$-equivalent to $\lambda
_{i}-\lambda _{j}$, $3\leq i<j$.

\item If $i>3$ then $\langle \gamma ,\lambda _{i}\pm \lambda _{j}\rangle
\neq 0$ but $\langle \gamma ,\lambda _{1}-\lambda _{2}\rangle =0$ where $%
\gamma =\lambda _{i-1}-\lambda _{i}$, and again $\lambda _{1}-\lambda _{2}$
is not $M$-equivalent to $\lambda _{3}-\lambda _{j}$, $3<j$.

\item If $i=3$ and $j=l+1$ then $\lambda _{4}-\lambda _{l+1}$ is a root
orthogonal to $\lambda _{1}-\lambda _{2}$ but not orthogonal to $\lambda
_{3}-\lambda _{l+1}$.
\end{enumerate}

Finally if $l=1$ there is just one positive root. If $l=2$ then the positive
roots are not orthogonal to each other so by Corollary \ref{cormequivalent}
they are not $M$-equivalent.

\subsubsection{$A_{3}$}

The $M$-equivalence classes on the positive roots are $\{\lambda
_{1}-\lambda _{2},\lambda _{3}-\lambda _{4}\}$, $\{\lambda _{1}-\lambda
_{3},\lambda _{2}-\lambda _{4}\}$ and $\{\lambda _{1}-\lambda _{4},\lambda
_{2}-\lambda _{3}\}$.

In this case the unique root orthogonal to $\alpha =\lambda _{1}-\lambda
_{2} $ is $\lambda _{3}-\lambda _{4}$ and hence, by Corollary \ref%
{cormequivalent}, $\lambda _{3}-\lambda _{4}$ is the only candidate to be $M$%
-equivalent to $\lambda _{1}-\lambda _{2}$. To see that indeed $\lambda
_{3}-\lambda _{4}\sim _{M}\lambda _{1}-\lambda _{2}$ note that a root $%
\gamma =\lambda _{i}-\lambda _{j}$ with $\left( i,j\right) \neq \left(
1,2\right) $ or $\left( 3,4\right) $ is not orthogonal to $\lambda
_{1}-\lambda _{2}$ neither to $\lambda _{3}-\lambda _{4}$. So that the
Killing numbers $\frac{2\langle \gamma ,\lambda _{1}-\lambda _{2}\rangle }{%
\langle \gamma ,\gamma \rangle }$ and $\frac{2\langle \gamma ,\lambda
_{3}-\lambda _{4}\rangle }{\langle \gamma ,\gamma \rangle }$ are $\pm 1$,
that is, \ the condition of Proposition \ref{propmequivalent} is satisfied
showing that $\lambda _{3}-\lambda _{4}\sim _{M}\lambda _{1}-\lambda _{2}$.
By applying the Weyl group (permutation group) we see that the classes of $M$%
-equivalences are as stated.

\subsection{Diagram $B_{l}$, $l\geq 2$}

In the standard realization of $B_{l}=\mathfrak{so}\left( l,l+1\right) $ the
positive roots are written as $\lambda _{i}\pm \lambda _{j}$, $1\leq i<j\leq
l$ and $\lambda _{i}$, $1\leq i\leq l$. These are the long and short roots
respectively.

The $M$-equivalence classes depend on the rank $l$, according to the
following cases:

\subsubsection{$B_{l}$, $l\geq 5$}

The $M$-equivalence classes on the positive roots are $\{\lambda
_{i}-\lambda _{j},\lambda _{i}+\lambda _{j}\}$ and $\{\lambda _{i}\}$, $%
1\leq i<j\leq l$.

We find the equivalence classes of the long and short roots.

Take long root $\lambda _{1}-\lambda _{2}$. We must check $M$-equivalence
only for the roots orthogonal to it, namely $\lambda _{1}+\lambda _{2}$, $%
\lambda _{i}\pm \lambda _{j}$ and $\lambda _{i}$ with $3\leq i<j$. The roots 
$\lambda _{i}$, $3\leq i$, are not $M$-equivalent to $\lambda _{1}-\lambda
_{2}$. In fact, $\lambda _{i}\pm \lambda _{i+1}$ is a root because $l\geq 5 $%
. Now, $\langle \lambda _{i},\lambda _{i}\pm \lambda _{i+1}\rangle \neq 0$
and the Killing number $\langle \left( \lambda _{i}\pm \lambda _{i+1}\right)
^{\vee },\lambda _{i}\rangle =\pm 1$ because $\lambda _{i}$ is a short root.
Since $\langle \lambda _{i},\lambda _{1}-\lambda _{2}\rangle =0$ the
condition of Proposition \ref{propmequivalent} is violated by $\gamma
=\lambda _{i}\pm \lambda _{i+1}$. The same argument used in the $A_{l}$ case
show that $\lambda _{i}\pm \lambda _{j}$, $3\leq i<j$, is not $M$-equivalent
to $\lambda _{1}-\lambda _{2}$ (when $l\geq 5$). On the other hand $\lambda
_{1}-\lambda _{2}\sim _{M}\lambda _{1}+\lambda _{2}$ because for any root $%
\gamma $ it holds $\langle \gamma ,\lambda _{1}-\lambda _{2}\rangle =\pm
\langle \gamma ,\lambda _{1}+\lambda _{2}\rangle $. It follows that $%
\{\lambda _{1}-\lambda _{2},\lambda _{1}+\lambda _{2}\}$ is an $M$%
-equivalence class. To conclude this case we note that $w\in \mathcal{W}$
acts on $\lambda _{i}$ by a permutation followed by a change of sign, that
is, $w\lambda _{i}=\pm \lambda _{j}$, for some index $j$. Hence $\lambda
_{i}-\lambda _{j}\sim _{M}\lambda _{i}+\lambda _{j}$, $1\leq i<j\leq l$ and
the sets $\{\lambda _{i}-\lambda _{j},\lambda _{i}+\lambda _{j}\}$ are the
only $M$-equivalence classes containing a long root.

By the previous paragraph no long root is $M$-equivalent to the short root $%
\lambda _{i}$. Finaly two short roots $\lambda _{i}$ and $\lambda _{j}$, $%
i\neq j$, are not $M$-equivalent. For example $\gamma =\lambda _{i}+\lambda
_{k}$, $k\neq i,j$ satisfies $\langle \gamma ^{\vee },\lambda _{i}\rangle =1 
$ while $\langle \gamma ^{\vee },\lambda _{j}\rangle =0$.

\subsubsection{$B_{4}$}

The $M$-equivalence classes on the positive roots are $\{\lambda
_{1}-\lambda _{2},\lambda _{1}+\lambda _{2},\lambda _{3}-\lambda
_{4},\lambda _{3}+\lambda _{4}\}$, $\{\lambda _{1}-\lambda _{3},\lambda
_{1}+\lambda _{3},\lambda _{2}-\lambda _{4},\lambda _{2}+\lambda _{4}\}$, $%
\{\lambda _{1}-\lambda _{4},\lambda _{1}+\lambda _{4},\lambda _{2}-\lambda
_{3},\lambda _{2}+\lambda _{3}\}$, and the short roots $\{\lambda _{i}\}$, $%
1\leq i\leq 4$.

The difference from the general case is that $\lambda _{3}-\lambda _{4}\sim
_{M}\lambda _{1}-\lambda _{2}$. In fact if $\lambda _{i}$ is a short root
then $\langle \lambda _{i}{}^{\vee },\lambda _{1}-\lambda _{2}\rangle $ and $%
\langle \lambda _{i}{}^{\vee },\lambda _{3}-\lambda _{4}\rangle $ equals $0$
or $2$. Also if $\gamma $ is a long root different from $\lambda
_{1}-\lambda _{2}$ and $\lambda _{3}-\lambda _{4}$ then $\langle \gamma
^{\vee },\lambda _{1}-\lambda _{2}\rangle $ and $\langle \gamma ^{\vee
},\lambda _{3}-\lambda _{4}\rangle $ equals $\pm 1$.

Again the same arguments show that a long root and a short root as well as
two short roots are not $M$-equivalent.

\subsubsection{$B_{3}$}

The $M$-equivalence classes on the positive roots are $\{\lambda
_{1}-\lambda _{2},\lambda _{1}+\lambda _{2},\lambda _{3}\}$, $\{\lambda
_{1}-\lambda _{3},\lambda _{1}+\lambda _{3},\lambda _{2}\}$ and $\{\lambda
_{2}-\lambda _{3},\lambda _{2}+\lambda _{3},\lambda _{1}\}$.

Here $\lambda _{3}\sim _{M}\lambda _{1}-\lambda _{2}$. The point is that if $%
\gamma \neq \lambda _{1}-\lambda _{2}$ is a long root then $\langle \gamma
^{\vee },\lambda _{1}-\lambda _{2}\rangle =\pm 1$ and since $\gamma $ cannot
be orthogonal to $\lambda _{3}$ we have $\langle \gamma ^{\vee },\lambda
_{3}\rangle =\pm 1$ as well. On the other hand if $\gamma $ is short then $%
\langle \gamma ^{\vee },\lambda _{1}-\lambda _{2}\rangle $ and $\langle
\gamma ^{\vee },\lambda _{3}\rangle $ are even.

\subsubsection{$B_{2}$}

The $M$-equivalence classes on the positive roots are the long roots $%
\{\lambda _{1}-\lambda _{2},\lambda _{1}+\lambda _{2}\}$ and the short roots 
$\{\lambda _{1},\lambda _{2}\}$.

\subsection{Diagram $C_{l}$, $l\geq 3$}

In the standard realization of $C_{l}=\mathfrak{sp}\left( l,\mathbb{R}%
\right) $ the positive roots are written as $\lambda _{i}\pm \lambda _{j}$, $%
1\leq i<j\leq l$ and $2\lambda _{i}$, $1\leq i\leq l$. These are the short
and long roots respectively.

Here any two long roots $2\lambda _{i}$ and $2\lambda _{j}$ are $M$%
-equivalent. In fact, for any root $\gamma $ the Killing number $\langle
\gamma ^{\vee },2\lambda _{i}\rangle $ is even ($0$ or $\pm 2$). In fact, if 
$\gamma $ is a short root then $\langle \gamma ^{\vee },2\lambda _{i}\rangle 
$ is either $0$ (orthogonal roots) or $\pm 2$ (Killing number between a
short root and a long root). On the other hand two long roots are either
equal or orthogonal.

As in the previous diagrams the $M$-equivalence classes increase for small
ranks. For $C_{l}$ the exception is when $l=4$.

\subsubsection{$C_{l}$, $l\neq 4$}

The $M$-equivalence classes are $\{\lambda _{i}-\lambda _{j},\lambda
_{i}+\lambda _{j}\}$ and the set of long roots $\{2\lambda _{1},\ldots
,2\lambda _{l}\}$.

The roots orthogonal to the short root $\lambda _{1}-\lambda _{2}$ are $%
\lambda _{1}+\lambda _{2}$, $\lambda _{i}\pm \lambda _{j}$ and $2\lambda
_{i} $ with $3\leq i<j$. As in the $B_{l}$ case (with $l\geq 5$) the roots $%
\lambda _{i}\pm \lambda _{j}$, $3\leq i<j$, are not $M$-equivalent to $%
\lambda _{1}-\lambda _{2}$. On the other hand if $3\leq i$ then $\gamma
=\lambda _{1}-\lambda _{j}$ with $j\neq i$ violates the criterion of
Proposition \ref{propmequivalent} for $M$-equivalence between $\lambda
_{1}-\lambda _{2}$ and $2\lambda _{i}$. In fact, $\langle \gamma ^{\vee
},\lambda _{1}-\lambda _{2}\rangle =\pm 1$ (non orthogonal roots of same
length) and $\langle \gamma ^{\vee },2\lambda _{i}\rangle =0$.

Since the long roots are equivalente to each other it follows that $\lambda
_{1}-\lambda _{2}$ is not $M$-equivalent to any long root. Hence we get the
classes stated above.

These arguments remain true if $l=3$. (Diferently from $B_{3}$ in $C_{3}$
long roots are not $M$-equivalent to short roots.)

\subsubsection{$C_{4}$}

The $M$-equivalence classes are $\{\lambda _{1}-\lambda _{2},\lambda
_{1}+\lambda _{2},\lambda _{3}-\lambda _{4},\lambda _{3}+\lambda _{4}\}$, $%
\{\lambda _{1}-\lambda _{3},\lambda _{1}+\lambda _{3},\lambda _{2}-\lambda
_{4},\lambda _{2}+\lambda _{4}\}$, $\{\lambda _{1}-\lambda _{4},\lambda
_{1}+\lambda _{4},\lambda _{2}-\lambda _{3},\lambda _{2}+\lambda _{3}\}$ and
the long roots $\{2\lambda _{1},2\lambda _{2},2\lambda _{3},2\lambda _{4}\}$.

This is seen as in $B_{4}$ where $\lambda _{3}-\lambda _{4}\sim _{M}\lambda
_{1}-\lambda _{2}$.

\subsection{Diagram $D_{l}$, $l\geq 4$}

In the standard realization of $D_{l}=\mathfrak{so}\left( l,l\right) $ the
positive roots are written as $\lambda _{i}\pm \lambda _{j}$, $1\leq i<j\leq
l$.

\subsubsection{$D_{l}$, $l>4$}

The $M$-equivalence classes on the positive roots are $\{\lambda
_{i}-\lambda _{j},\lambda _{i}+\lambda _{j}\}$, $1\leq i<j\leq l$.

This is verified by arguments similar to the $B_{l}$ case, simplified by the
fact that the roots have the same length.

First the only root $M$-equivalent to $\lambda _{1}-\lambda _{2}$ is $%
\lambda _{1}+\lambda _{2}$. In fact, the roots orthogonal to $\lambda
_{1}-\lambda _{2}$ are $\lambda _{1}+\lambda _{2}$ and $\lambda _{i}\pm
\lambda _{j}$, $3\leq i<j$. A root $\lambda _{i}\pm \lambda _{j}$ with $%
3\leq i<j$ is not $M$-equivalent to $\lambda _{1}-\lambda _{2}$ by the
following reasons:

\begin{enumerate}
\item If $j<l$ and $\gamma =\lambda _{j}-\lambda _{j+1}$ then $\langle
\gamma ,\lambda _{1}-\lambda _{2}\rangle =0$ and $\langle \gamma ,\lambda
_{i}\pm \lambda _{j}\rangle \neq 0$ which implies that $\langle \gamma
^{\vee },\lambda _{i}\pm \lambda _{j}\rangle =\pm 1$. Thus by Proposition %
\ref{propmequivalent} $\lambda _{1}-\lambda _{2}$ is not $M$-equivalent to $%
\lambda _{i}\pm \lambda _{j}$.

\item If $i>3$ and $\gamma =\lambda _{i-1}-\lambda _{i}$ then $\langle
\gamma ,\lambda _{1}-\lambda _{2}\rangle =0$ and $\langle \gamma ^{\vee
},\lambda _{i}\pm \lambda _{j}\rangle =\pm 1$.

\item Since $l>4$, $\lambda _{4}-\lambda _{l}$ is a root satisfying $\langle
\gamma ,\lambda _{1}-\lambda _{2}\rangle =0$ and $\langle \gamma ^{\vee
},\lambda _{3}\pm \lambda _{l}\rangle =\pm 1$.
\end{enumerate}

Finally $\lambda _{1}-\lambda _{2}\sim _{M}\lambda _{1}+\lambda _{2}$,
because $\langle \gamma ,\lambda _{1}-\lambda _{2}\rangle =0$ if and only if 
$\langle \gamma ,\lambda _{1}+\lambda _{2}\rangle =0$ for any root $\gamma $%
. Also, if $\gamma $ is not orthogonal to both roots then the Killing
numbers are $\pm 1$, since the roots have the same length.

Since the Weyl group is transitive on the set of roots we get the
equivalence classes stated above.

\subsubsection{$D_{4}$}

The $M$-equivalence classes on the positive roots are $\{\lambda
_{1}-\lambda _{2},\lambda _{1}+\lambda _{2},\lambda _{3}-\lambda
_{4},\lambda _{3}+\lambda _{4}\}$, $\{\lambda _{1}-\lambda _{3},\lambda
_{1}+\lambda _{3},\lambda _{2}-\lambda _{4},\lambda _{2}+\lambda _{4}\}$ and 
$\{\lambda _{1}-\lambda _{4},\lambda _{1}+\lambda _{4},\lambda _{2}-\lambda
_{3},\lambda _{2}+\lambda _{3}\}$.

In this case, appart from $\lambda _{1}+\lambda _{2}$ the roots $\lambda
_{3}-\lambda _{4}$ and $\lambda _{3}+\lambda _{4}$ are $M$-equivalent to $%
\lambda _{1}-\lambda _{2}$ (see the discussion for $B_{4}$). Hence an
application of the Weyl group yield the stated classes.

\subsection{Diagram $G_{2}$}

The $M$-equivalence classes on the positive roots are the pairs $\{\alpha
_{1},\alpha _{1}+2\alpha _{2}\}$, $\{\alpha _{1}+\alpha _{2},\alpha
_{1}+3\alpha _{2}\}$ and $\{\alpha _{2},2\alpha _{1}+3\alpha _{2}\}$ where $%
\alpha _{1}$ and $\alpha _{2}$ are the simple roots with $\alpha _{1}$ the
long one.

The reason is that these are the only pairs of positive roots orthogonal to
each other. Moreover if two roots are not orthogonal then their Killing are
odd ($\pm 1$ ou $\pm 3$).

\subsection{Diagrams $E_{6}$, $E_{7}$ and $E_{8}$}

For these diagrams the $M$-equivalence classes on the positive roots are
singletons.

Since these diagrams are simply-laced it is enough to find a positive root
which is not $M$-equivalent to any other positive root.

In any of the diagrams $E_{6}$, $E_{7}$ and $E_{8}$ we choose the highest
root $\mu $. To check that $\{\mu \}$ is an $M$-equivalence class we prove
the

\begin{itemize}
\item \textbf{Claim:} For every $\beta >0$ with $\langle \mu ,\beta \rangle
=0$ there exists $\gamma \neq \beta $ such that $\langle \mu ,\gamma \rangle
=0$ and $\langle \beta ,\gamma \rangle \neq 0$.
\end{itemize}

From the claim we get $\langle \gamma ^{\vee },\mu \rangle =0$ and $\langle
\gamma ^{\vee },\beta \rangle $ odd because the diagrams are simply laced.
Hence, by Proposition \ref{propmequivalent}, no $\beta $ orthogonal to $\mu $
is $M$-equivalent to $\mu $. By Corollary \ref{cormequivalent} we conclude
that $\{\mu \}$ is an $M$-equivalence class.

Now the roots orthogonal to the highest root $\mu $ have the following
simple description: Denote by $\Sigma =\{\alpha _{1},\ldots ,\alpha _{l}\}$
the simple system of roots, and let $\{\omega _{1},\ldots ,\omega _{l}\}$ be
the fundamental weights, defined by 
\begin{equation*}
\langle \alpha _{i}^{\vee },\omega _{j}\rangle =\delta _{ij}.
\end{equation*}%
It is known that in the diagrams $E_{6}$, $E_{7}$ and $E_{8}$ the highest
root $\mu =\omega _{i}$ for some fundamental weight. (The formula for $\mu $
in terms of the fundamental weights can be read off from the affine Dynkin
diagrams. The extra root is precisely $-\mu $, see \cite{he}, Chapter X,
Table of Diagrams S(A).) Let $\alpha =b_{1}\alpha _{1}+\cdots +b_{l}\alpha
_{l}$, $b_{i}\geq 0$, be a positive root. Since $\mu =\omega _{i}$ we have
by definition 
\begin{equation*}
\langle \alpha ,\mu \rangle =\frac{\langle \alpha _{i},\alpha _{i}\rangle }{2%
}a_{i}b_{i}.
\end{equation*}%
So that $\langle \alpha ,\mu \rangle =0$ if and only if $a_{i}b_{i}=0$.
Therefore the roots orthogonal to $\mu $ are those spanned by $\Sigma
\setminus \{\alpha _{i}\}$. This set of roots is a root system whose Dynkin
diagram is the subdiagram of $\Sigma $ given by $\Sigma \setminus \{\alpha
_{i}\}$. A glance at the affine Dynkin diagrams provides the diagrams $%
\Sigma \setminus \{\alpha _{i}\}$, namely,

\begin{itemize}
\item $\Sigma \setminus \{\alpha _{i}\}=A_{5}$ if $\Sigma =E_{6}$.

\item $\Sigma \setminus \{\alpha _{i}\}=D_{6}$ if $\Sigma =E_{7}$.

\item $\Sigma \setminus \{\alpha _{i}\}=E_{7}$ if $\Sigma =E_{8}$.
\end{itemize}

Now it is clear that in any of the root systems spanned by $\Sigma \setminus
\{\alpha _{i}\}$ ($A_{5}$, $D_{6}$ or $E_{7}$), the conclusion of the claim
holds, that is, given $\beta $ there exists $\gamma $ with $\langle \beta
,\gamma \rangle \neq 0$. This concludes the proof that the $M$-equivalence
classes on the positive roots are singletons.

\subsection{$F_{4}$}

The $24$ positive roots of

%
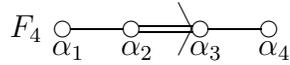
\begin{picture}(100,60)(0,0)

\put(0,27){$F_4$}
\put(20,30){\circle{6}}
\put(16,20){${\alpha}_1$}
\put(23,30){\line(1,0){20}}
\put(46,30){\circle{6}}
\put(42,20){${\alpha}_2$}
\put(49,31.2){\line(1,0){20}}
\put(49,28.8){\line(1,0){20}}
\put(72,30){\circle{6}}
\put(68,20){${\alpha}_3$}
\put(69,30){\line(-1,2){5}}
\put(69,30){\line(-1,-2){5}}
\put(75,30){\line(1,0){20}}
\put(98,30){\circle{6}}
\put(94,20){${\alpha}_4$}
\end{picture}%

\noindent%
split into the following $M$-equivalence classes:

\begin{itemize}
\item $12$ singletons $\{\alpha \}$ with $\alpha $ running through the set
of short roots.

\item $3$ sets of long roots $\{2\alpha _{1}+3\alpha _{2}+4\alpha
_{3}+2\alpha _{4},\alpha _{2},\alpha _{2}+2\alpha _{3},\alpha _{2}+2\alpha
_{3}+2\alpha _{4}\}$, $\{\alpha _{1}+3\alpha _{2}+4\alpha _{3}+2\alpha
_{4},\alpha _{1}+\alpha _{2},\alpha _{1}+\alpha _{2}+2\alpha _{3},\alpha
_{1}+\alpha _{2}+2\alpha _{3}+2\alpha _{4}\}$ and

$\{\alpha _{1}+2\alpha _{2}+4\alpha _{3}+2\alpha _{4},\alpha _{1},\alpha
_{1}+2\alpha _{2}+2\alpha _{3},\alpha _{1}+2\alpha _{2}+2\alpha _{3}+2\alpha
_{4}\}$.
\end{itemize}

We let $\{\omega _{1},\omega _{2},\omega _{3},\omega _{4}\}$ be the
fundamental weights.

The fundamental weight $\omega _{4}$ is also the short\ positive root%
\begin{equation*}
\omega _{4}=\alpha _{1}+2\alpha _{2}+3\alpha _{3}+2\alpha _{4}.
\end{equation*}%
We look at its $M$-equivalence class by the same method of the $E_{l}$'s.
The set of roots orthogonal to the fundamental weight $\omega _{4}$ is
spanned by $\{\alpha _{1},\alpha _{2},\alpha _{3}\}$ which is a $B_{3}$
Dynkin diagram. Now if $\beta $ is a root of $B_{3}$ then there exists a
root $\gamma $ (in $B_{3}$) such that $\langle \gamma ^{\vee },\beta \rangle 
$ is odd. It follows by Proposition \ref{propmequivalent} and its Corollary %
\ref{cormequivalent} that $\{\omega _{4}\}$ is an $M$-equivalence class.
This gives the classes of the short roots.

As to the long roots we first recall that they form a $D_{4}$ root system
(see ????). Now if $\gamma $ is a short root and $\alpha $ a long root then $%
\langle \gamma ^{\vee },\alpha \rangle $ is even. Hence to check if $\alpha
\sim _{M}\beta $ for the two long roots $\alpha $ and $\beta $ it is enough
to test the condition of Proposition \ref{propmequivalent} when $\gamma $ is
also a long root. This means that two long roots are $M$-equivalent if and
only if they are equivalent as roots of $D_{4}$. Since no short root is $M$%
-equivalent to a long root we conclude that the classes of $D_{4}$ are also $%
M$-equivalence classes in $F_{4}$. These are the three sets with four
orthogonal roots each as stated. (To get these sets start with the highest
root $\omega _{1}=2\alpha _{1}+3\alpha _{2}+4\alpha _{3}+2\alpha _{4}$. Then
the first set is $\omega _{1}$ together with the long roots orthogonal to
it. The next two sets are obtained by applying first the reflection $%
r_{\alpha _{1}}$ and then $r_{\alpha _{2}}$.)

\ \ 

\section{Auxiliary lemmas\label{seclemaux}}

In this section we prove some lemmas to be used later in the determination
of the irreducible $K_{\Theta }$-invariant subspaces of $\mathfrak{n}%
_{\Theta }^{-}$. We choose once and for all a generator $E_{\alpha }\in 
\mathfrak{g}_{\alpha }$ for each root space.

Recall that in Section \ref{secisotro} we denoted the irreducible components
for the adjoint representation of $\mathfrak{z}_{\Theta }$ on $\mathfrak{n}%
_{\Theta }^{-}$ by $V_{\Theta }^{\sigma }$. Write $\Pi _{\Theta }^{\sigma
}\subset \Pi ^{-}\backslash \langle \Theta \rangle ^{-}$ for the set of
roots such that 
\begin{equation*}
V_{\Theta }^{\sigma }=\sum_{\alpha \in \Pi _{\Theta }^{\sigma }}\mathfrak{g}%
_{\alpha }.
\end{equation*}%
The subgroup $K_{\Theta }$ leave invariant $V_{\Theta }^{\sigma }$, hence $%
\Pi _{\Theta }^{\sigma }$ is $\mathcal{W}_{\Theta }$-invariant. \ 

Our first results give conditions ensuring that a \ $\mathfrak{z}_{\Theta }$%
-invariant subspace is $K_{\Theta }$-irreducible.

\begin{lema}
\label{leminterinvMclass}Let $V=\sum_{\alpha \in \Pi _{V}}\mathfrak{g}%
_{\alpha }$ be a $\mathfrak{z}_{\Theta }$-invariant subspace and suppose
that $W\subset V$ is a $K_{\Theta }$-invariant subspace. Take $%
X=\sum_{\alpha \in \Pi _{V}}a_{\alpha }E_{\alpha }\in W$ and let $\alpha $
be a root such that $a_{\alpha }\neq 0$. Define 
\begin{equation*}
V_{[\alpha ]_{M}}=\sum_{\beta \in \left[ \alpha \right] _{M}\cap \Pi _{V}}%
\mathfrak{g}_{\beta }.
\end{equation*}

Then $W\cap V_{[\alpha ]_{M}}\neq \{0\}$.
\end{lema}

\begin{profe}
Let $c_{X,\alpha }$ be the cardinality of $\{\beta \notin \lbrack \alpha
]_{M}:a_{\beta }\neq 0\}$. If $c_{X,\alpha }=0$ we are done. Otherwise we
find $0\neq Y\in U$ with $c_{Y,\alpha }<c_{X,\alpha }$. In fact, if $%
c_{X,\alpha }>0$ then there are $\beta \notin \lbrack \alpha ]_{M}$ with $%
a_{\alpha },a_{\beta }\neq 0$. So that there exists $m\in M$ with $%
mE_{\alpha }=E_{\alpha }$ and $mE_{\beta }=-E_{\beta }$.

Now, $M\subset K_{\Theta }$, hence $Y=X+mX\in W$. Clearly $Y\neq 0$ and
since the $\beta $ component of $Y$ is zero we have $c_{Y,\alpha
}<c_{X,\alpha }$. Repeating this argument successively we arrive at $Z\in W$
such that $c_{Z,\alpha }=0$, concluding the proof.
\end{profe}

\begin{lema}
\label{lemirreducible}Take a subset $\Pi _{V}\subset \Pi ^{-}\backslash
\langle \Theta \rangle ^{-}$ such that the subspace $V=\sum_{\alpha \in \Pi
_{V}}\mathfrak{g}_{\alpha }$ is $\mathfrak{z}_{\Theta }$-invariant. Suppose
that

\begin{enumerate}
\item $\mathcal{W}_{\Theta }$ acts transitively on $\Pi _{V}$, and

\item two different roots in $\Pi _{V}$ are not $M$-equivalent.
\end{enumerate}

Then $V$ is $K_{\Theta }$-irreducible.
\end{lema}

\begin{profe}
Let $U\subset V$ be a nontrivial $K_{\Theta }$-invariant subspace. By
transitivity of $\mathcal{W}_{\Theta }\subset K_{\Theta }$ it is enough to
prove that $U$ contains a root space $\mathfrak{g}_{\alpha }$, $\alpha \in
\Pi _{V}$. But this follows by the previous lemma and the assumption that
different roots in $\Pi _{V}$ are not $M$-equivalent.
\end{profe}

As a complement of the above lemma we exhibit next general cases where $%
\mathcal{W}_{\Theta }$ acts transitively on sets of roots.

\begin{lema}
\label{lemtransimplelong}Let $\Pi _{\Theta }^{\sigma }\subset \Pi
^{-}\backslash \langle \Theta \rangle ^{-}$ be the set of roots
corresponding to an irreducible component $V_{\Theta }^{\sigma
}=\sum_{\alpha \in \Pi _{\Theta }^{\sigma }}\mathfrak{g}_{\alpha }$. In each
of the following cases $\mathcal{W}_{\Theta }$ acts transitively on $\Pi
_{\Theta }^{\sigma }$.

\begin{enumerate}
\item The Dynkin diagram of $\mathfrak{g}$ has only simple edges ($A_{l}$, $%
D_{l}$, $E_{6}$, $E_{7}$ and $E_{8}$).

\item For the diagrams $B_{l}$, $C_{l}$ and $F_{4}$ there are the cases:

\begin{enumerate}
\item The roots in $\Theta \subset \Sigma $ are long.

\item The roots in $\Pi _{\Theta }^{\sigma }$ are short.
\end{enumerate}
\end{enumerate}
\end{lema}

\begin{profe}
Let $\mu $ be the highest root of $\Pi _{\Theta }^{\sigma }$. By
representation theory we know that any other root $\beta \in \Pi _{\Theta
}^{\sigma }$ (weight of the representation) is given by 
\begin{equation*}
\beta =\mu -\alpha _{1}-\cdots -\alpha _{k}
\end{equation*}%
with $\alpha _{i}\in \Theta $ and such that any partial difference $\mu
-\alpha _{1}-\cdots -\alpha _{i}$, $i\leq k$, is also a root. Fix $i\leq k$
put $\delta =\mu -\alpha _{1}-\cdots -\alpha _{i-1}$ and let $r_{\alpha
_{i}} $ be the reflection with respect to $\alpha _{i}$.

We claim that $\delta -\alpha _{i}=r_{\alpha _{i}}\left( \delta \right) $.
This follows by the Killing formula applied to the string of roots $\delta
+k\alpha _{i}$. There are the following cases:

\begin{enumerate}
\item In the simply laced diagrams of (1) the Killing number 
\begin{equation*}
\langle \alpha _{i}^{\vee },\delta \rangle =\frac{2\langle \alpha
_{i},\delta \rangle }{\langle \alpha _{i},\alpha _{i}\rangle }
\end{equation*}%
is $0$ or $\pm 1$. Since $\delta -\alpha _{i}$ is a root we have $\langle
\alpha _{i}^{\vee },\delta \rangle =1$, and hence $r_{\alpha _{i}}\left(
\delta \right) =\delta -\alpha _{i}$.

\item If the roots in $\Theta $ are long as in (2a) then $\alpha _{i}$ is a
long root implying that $\langle \alpha _{i}^{\vee },\delta \rangle $ is $0$
or $\pm 1$. Again, the fact that $\delta -\alpha _{i}$ is a root implies
that $\langle \alpha _{i}^{\vee },\delta \rangle =1$, so that $r_{\alpha
_{i}}\left( \delta \right) =\delta -\alpha _{i}$.

\item If the roots in $\Pi _{\Theta }^{\sigma }$ are short in a double laced
diagram as in (2b) then $\delta $ and $\delta -\alpha _{i}$ are a short
roots. If $\alpha _{i}$ is a long root then $\delta $ and $\delta -\alpha
_{i}$ are the only roots of the form $\delta +k\alpha _{i}$, $k\in \mathbb{Z}
$. Hence by the Killing formula $\langle \alpha _{i}^{\vee },\delta \rangle
=1$, that is $r_{\alpha _{i}}\left( \delta \right) =\delta -\alpha _{i}$.

On the other hand if $\alpha _{i}$ is short then there are two possibilities
for the string of roots $\delta +k\alpha _{i}$: i) $\delta -\alpha _{i}$, $%
\delta $ and $\delta +\alpha _{i}$ are roots in which case $\langle \alpha
_{i},\delta \rangle =0$ and $\delta -\alpha _{i}$ and $\delta +\alpha _{i}$
are long roots; ii) $\delta -\alpha _{i}$ and $\delta $ are roots and $%
\langle \alpha _{i},\delta \rangle =1$. The first case is ruled out because
otherwise we would have the long roots $\delta -\alpha _{i},\delta +\alpha
_{i}\in \Pi _{\Theta }^{\sigma }$, contradicting the assumption. Therefore $%
\langle \alpha _{i},\delta \rangle =1$, that is, $r_{\alpha _{i}}\left(
\delta \right) =\delta -\alpha _{i}$.
\end{enumerate}

Since $r_{\alpha _{i}}\in \mathcal{W}_{\Theta }$, it follows by induction
that $\beta $ belongs to the $\mathcal{W}_{\Theta }$-oribt of $\mu $,
proving transitivity of $\mathcal{W}_{\Theta }$. \ 
\end{profe}

We turn now to the equivalence of irreducible representations.

\begin{lema}
\label{lemnotequivalent}Let $V_{\Theta }^{\sigma }$ and $V_{\Theta }^{\tau }$
be $\mathfrak{z}_{\Theta }$-irreducible components. Suppose that there
exists $\alpha \in \Pi _{\Theta }^{\sigma }$ which is not $M$-equivalent to
any $\beta \in \Pi _{\Theta }^{\tau }$. Then $V_{\Theta }^{\sigma }$ and $%
V_{\Theta }^{\tau }$ are not $K_{\Theta }$-equivalent.
\end{lema}

\begin{profe}
Suppose to the contrary that there exists an isomorphism $T:V_{\Theta
}^{\sigma }\rightarrow V_{\Theta }^{\tau }$ intertwining the $K_{\Theta }$%
-representations. In particular 
\begin{equation*}
TmX=mTX,
\end{equation*}%
for all $m\in M\subset K_{\Theta }$ and $X\in V_{\Theta }^{\sigma }$.

Take $0\neq E_{\alpha }\in \mathfrak{g}_{\alpha }$. Then for every $m\in M$
we have $mE_{\alpha }=\allowbreak \varepsilon _{m}E_{\alpha }$ with $%
\varepsilon _{m}=\pm 1$. Write 
\begin{equation*}
TE_{\alpha }=\sum_{\beta \in \Pi _{\Theta }^{\tau }}a_{\beta }E_{\beta }.
\end{equation*}%
Then for $m\in M$ we have 
\begin{equation*}
\varepsilon _{m}\sum_{\beta \in \Pi _{\Theta }^{\tau }}a_{\beta }E_{\beta
}=\varepsilon _{m}TE_{\alpha }=mTE_{\alpha }=\sum_{\beta \in \Pi _{\Theta
}^{\tau }}a_{\beta }mE_{\beta }.
\end{equation*}%
Since $mE_{\beta }=\pm E_{\beta }$ and the set $E_{\beta }$ is linearly
independent, it follows that $mE_{\beta }=\varepsilon _{m}E_{\beta }$ if $%
a_{\beta }\neq 0$. For any such $\beta $ the representation of $M$ on $%
\mathfrak{g}_{\beta }$ is equivalent to the representation on $\mathfrak{g}%
_{\alpha }$. This contradicts the assumption that $\alpha $ is not $M$%
-equivalent to $\beta \in \Pi _{\Theta }^{\tau }$.
\end{profe}

The next statement gives a sufficient condition for equivalence.

\begin{proposicao}
\label{propequivalent}Let $V_{\Theta }^{\sigma }$ and $V_{\Theta }^{\tau }$
be $\mathfrak{z}_{\Theta }$-irreducible components. Suppose that there is a
bijection $\iota :\Pi _{\Theta }^{\sigma }\rightarrow \Pi _{\Theta }^{\tau }$
such that $\mathfrak{g}_{\alpha }$ and $\mathfrak{g}_{\iota \left( \alpha
\right) }$ are $M$-equivalent for every $\alpha \in \Pi _{\Theta }^{\sigma }$%
. Assume also that the linear map $T:V_{\Theta }^{\sigma }\rightarrow
V_{\Theta }^{\tau }$, given by $TE_{\alpha }=E_{\iota \left( \alpha \right)
} $, commutes with $\mathrm{ad}\left( X\right) $, $X\in \mathfrak{k}_{\Theta
}$ for every $\alpha \in \Pi _{\Theta }^{\sigma }$. Then $T$ is an
intertwining operator for the $K_{\Theta }$-representations on $V_{\Theta
}^{\sigma }$ and $V_{\Theta }^{\tau }$. Moreover the subspaces 
\begin{equation*}
V_{[(x,y)]}=\left\{ xX+yTX:X\in V_{\Theta }^{\sigma }\right\} ,
\end{equation*}%
where $[(x,y)]\in \mathbb{RP}^{2}$, are the only $K_{\Theta }$-invariant
subspaces in $V_{\Theta }^{\sigma }\oplus V_{\Theta }^{\tau }$.
\end{proposicao}

\begin{profe}
The first assumption implies that $T$ intertwines the $M$-representations,
while the second assumption means that $T$ intertwines the representations
of $\left( K_{\Theta }\right) _{0}$. Since $K_{\Theta }=M\left( K_{\Theta
}\right) _{0}$, we conclude that $T$ is in fact an intertwining operator for
the $K_{\Theta }$-representations.

This implies that $V_{[(x,y)]}$ is a $K_{\Theta }$-invariant subspace in $%
V_{\Theta }^{\sigma }\oplus V_{\Theta }^{\tau }$ for any $[(x,y)]\in \mathbb{%
RP}^{2}$.

Now, if $V$ is a $K_{\Theta }$-invariant subspace in $V_{\Theta }^{\lambda
}\oplus V_{\Theta }^{\mu }$ different from $V_{\Theta }^{\sigma
}=V_{[(1,0)]} $ or $V_{\Theta }^{\tau }=V_{[(0,1)]}$, then there exist a
linear isomorphism $L:V_{\Theta }^{\sigma }\rightarrow V_{\Theta }^{\tau }$
such that 
\begin{equation*}
V=\{X+LX:X\in V_{\Theta }^{\lambda }\}
\end{equation*}%
and 
\begin{equation*}
LkX=kLX,
\end{equation*}%
for every $k\in K_{\Theta }$. Since $M$ is a subset of $K_{\Theta }$ and
since $mE_{\alpha }=\varepsilon _{m}E_{\alpha }$, we can argue as in the
proof of the previous Lemma to show that 
\begin{equation*}
LE_{\alpha }=y_{i}E_{\iota \left( \alpha \right) }.
\end{equation*}%
Since $\mathcal{W}_{\Theta }$ acts transitively in the set of the directions 
$\{\mathfrak{g}_{1}^{\lambda },\ldots ,\mathfrak{g}_{n_{\lambda }}^{\lambda
}\}$, we conclude that $y_{i}=y$, is independent of the index $i\in
\{1,\ldots ,n_{\lambda }\}$. Thus we have that $V=V_{[(1,y)]}$, concluding
the proof.
\end{profe}

The previous results are complemented by the following standard basic fact
in representation theory.

\begin{proposicao}
\label{propsumirredcomp}Let $V$ be the space of a finite dimensional
representation of a group $L$. Suppose that 
\begin{equation*}
V=V_{1}\oplus \cdots \oplus V_{s}
\end{equation*}%
with $V_{i}$ invariant and irreducible. If the representations of $L$ on
different components $V_{i}$, $V_{j}$, $i\neq j$, are not equivalent then
the only $L$-invariant subspaces are sums of the components.
\end{proposicao}

\begin{profe}
(Sketch) If $\{0\}\neq W\subset V$ is an invariant subspace then the
projection $W_{i}$ to $V_{i}$ is invariant and hence either $\{0\}$ or $%
V_{i} $. Suppose that there are $i\neq j$ such that $W_{i}=V_{i}$ and $%
W_{j}=V_{j}$ and write $W_{ij}$ for the projection on $V_{i}\oplus V_{j}$.
Then $W_{ij}\cap V_{i}$ is $\{0\}$ or $V_{i}$. If $W_{ij}\cap V_{i}=V_{i}$
then $W_{ij}=V_{i}\oplus V_{j}$, which implies that $V_{i}\oplus
V_{j}\subset W$. Otherwise $W_{ij}\cap V_{i}=W_{ij}\cap V_{j}=\{0\}$. In
this case $W_{ij}$ is the graph of an isomorphism $V_{i}\rightarrow V_{j}$,
intertwining the representations on $V_{i}$ and $V_{j}$.
\end{profe}

Finally for several split simple Lie algebras the compact subalgebra $%
\mathfrak{k}$ is not simple. Via the next lemma we exploit this fact to get $%
K_{\Theta }$-invariant subspaces in $\mathfrak{n}_{\Theta }^{-}$.

\begin{lema}
\label{lemorbitnormal}Let $U\subset K$ be a normal subgroup and denote by $%
V\subset \mathfrak{n}_{\Theta }^{-}$ the tangent space to the $U$-orbit $%
U\cdot b_{\Theta }$ through the origin. Then $V$ is $K_{\Theta }$-invariant.
\end{lema}

\begin{profe}
The orbit $U\cdot b_{\Theta }$ is invariant by $K_{\Theta }$. In fact, if $%
u\cdot b_{\Theta }\in U\cdot b_{\Theta }$ and $k\in K_{\Theta }$ then $%
kuk^{-1}\in U$ so that $ku\cdot b_{\Theta }=kuk^{-1}\cdot kb_{\Theta
}=kuk^{-1}\cdot b_{\Theta }$ belongs to $U\cdot b_{\Theta }$. Hence its
tangent space at $b_{\Theta }$ is invariant by the isotropy representation.
\end{profe}

\section{Irreducible $K_{\Theta }$-invariant subspaces\label{secirreduc}}

Inn this section we describe \ the previous results to each diagram.

\subsection{Flags of $A_{l}=\mathfrak{sl}\left( l+1,\mathbb{R}\right) $}

As checked in Section \ref{secmequiv} no two different negative roots of $%
A_{l}$ are $M$-equivalent if $l\neq 3$. On the other hand by Lemma \ref%
{lemtransimplelong}, on any flag manifold of $A_{l}$, the subgroup $\mathcal{%
W}_{\Theta }$ acts transitively on each set of roots $\Pi _{\Theta }^{\sigma
}$ corresponding to an irreducible representation of $\mathfrak{z}_{\Theta }$
on $V_{\Theta }^{\sigma }$. Therefore by Lemma \ref{lemirreducible} we
conclude that $K_{\Theta }$ is irreducible on each $V_{\Theta }^{\sigma }$.
Looking again the $M$-equivalence classes we see that two different
irreducible subspaces are not $K_{\Theta }$-equivalent. Hence we get the
following description of the $K_{\Theta }$-invariant irreducible subspaces
in a flag manifold of $A_{l}$.

\begin{proposicao}
For any flag manifold $\mathbb{F}_{\Theta }$ of $A_{l}$, $l\neq 3$, the $%
K_{\Theta }$-invariant irreducible subspaces are the irreducible components $%
V_{\Theta }^{\sigma }$ for the $\mathfrak{z}_{\Theta }$ representation. Two
such representations are not $K_{\Theta }$-equivalent.
\end{proposicao}

The irreducible components $V_{\Theta }^{\sigma }$ are easily described in
terms of the matrices in $\mathfrak{sl}(n,\mathbb{R})$, $n=l+1$. In fact,
let 
\begin{equation*}
H_{\Theta }=\mathrm{diag}\{a_{1},\ldots ,a_{n}\}\in \mathfrak{a}_{\Theta
}\cap \mathrm{cl}\mathfrak{a}^{+}\quad a_{1}\geq \cdots \geq a_{n}
\end{equation*}%
be characteristic for $\Theta $. The multiplicities of the eigenvalues of $%
H_{\Theta }$ determine the sizes of a block decomposition of the $n\times n$
matrices. With respect to this decomposition the matrices in $\mathfrak{z}%
_{\Theta }$ are block diagonal while a block outside the diagonal determines
a $\mathfrak{z}_{\Theta }$-irreducible component. These are also the $%
K_{\Theta }$-irreducible components.

Now we look at the case $l=3$. The matrix 
\begin{equation*}
\left( 
\begin{array}{llll}
\ast & a & b & c \\ 
a & \ast & c & b \\ 
b & c & \ast & a \\ 
c & b & a & \ast%
\end{array}%
\right)
\end{equation*}%
summarizes the $M$-equivalence classes of $\mathfrak{sl}\left( 4,\mathbb{R}%
\right) $, where root spaces represented by the same letter are $M$%
-equivalent (see Section \ref{secmequiv}).

This shows that in the maximal flag manifold $\mathbb{F}_{\Theta }$, $\Theta
=\emptyset $, the $K_{\Theta }=M$ invariant irreducible subspaces are the
one-dimensional subspaces of $\mathfrak{g}_{21}\oplus \mathfrak{g}_{43}$, $%
\mathfrak{g}_{31}\oplus \mathfrak{g}_{42}$ or $\mathfrak{g}_{32}\oplus 
\mathfrak{g}_{41}$. The irreducible subspaces in the other flag manifolds
are easily obtained from this $M$-equivalence.

We discuss further the instructive case when $\mathbb{F}_{\Theta }=\mathrm{Gr%
}_{2}\left( 4\right) $, the Grassmannian of two dimensional subspaces of $%
\mathbb{R}^{4}$. In this case $\mathfrak{n}_{\Theta }^{-}$ is the subalgebra
of matrices writen in $2\times 2$ blocks as 
\begin{equation*}
X=\left( 
\begin{array}{ll}
0 & 0 \\ 
B & 0%
\end{array}%
\right) .
\end{equation*}%
The representations of $\mathfrak{z}_{\Theta }$ and $\mathfrak{g}\left(
\Theta \right) $ on $\mathfrak{n}_{\Theta }^{-}$ are irreducible. Here $%
K_{\Theta }=\mathrm{SO}\left( 2\right) \times \mathrm{SO}\left( 2\right) $
whose representation on $\mathfrak{n}_{\Theta }^{-}$ decomposes into two $2$%
-dimensional irreducible subspaces. This is due to the fact that $\mathfrak{%
so}\left( 4\right) =\mathfrak{so}\left( 3\right) _{1}\oplus \mathfrak{so}%
\left( 3\right) _{2}$ is a sum of two copies of $\mathfrak{so}\left(
3\right) $. The matrices in these components have the form%
\begin{equation}
\mathfrak{so}\left( 3\right) _{1}:\left( 
\begin{array}{cc}
A & -B^{T} \\ 
B & A%
\end{array}%
\right) \qquad \mathfrak{so}\left( 3\right) _{2}:\left( 
\begin{array}{cc}
A & -B^{T} \\ 
B & -A%
\end{array}%
\right)  \label{fordecomso4so3so3}
\end{equation}%
with $A+A^{T}=0$ where $B$ is symmetric with $\mathrm{tr}B=0$ for $\mathfrak{%
so}\left( 3\right) _{1}$ while 
\begin{equation}
B=\left( 
\begin{array}{cc}
a & -b \\ 
b & a%
\end{array}%
\right)  \label{forBcomplex}
\end{equation}%
for $\mathfrak{so}\left( 3\right) _{2}$. Hence by Lemma \ref{lemorbitnormal}%
, the tangent spaces $V_{i}$ to orbits of $\mathrm{SO}\left( 3\right)
_{i}=\langle \exp \mathfrak{so}\left( 3\right) _{i}\rangle $, $i=1,2$, are $%
K_{\Theta }$-invariant. The subspace $V_{i}$, $i=1,2$, is given by the
matrices in $\mathfrak{n}_{\Theta }^{-}$ with $B$ as $\mathfrak{so}\left(
3\right) _{1}$ or $\mathfrak{so}\left( 3\right) _{2}$, respectively.

\subsection{Flags of $B_{l}=\mathfrak{sl}\left( l+1,l\right) $}

In the standard realization $\mathfrak{sl}\left( l+1,l\right) $ is the
algebra of matrices%
\begin{equation*}
\left( 
\begin{array}{ccc}
0 & a & b \\ 
-b^{T} & A & B \\ 
-a^{T} & C & -A^{T}%
\end{array}%
\right) \qquad B+B^{T}=C+C^{T}=0.
\end{equation*}%
In this case $\mathfrak{a}$ is the subalgebra of matrices 
\begin{equation*}
\left( 
\begin{array}{ccc}
0 & 0 & 0 \\ 
0 & \Lambda & 0 \\ 
0 & 0 & -\Lambda%
\end{array}%
\right)
\end{equation*}%
with $\Lambda =\mathrm{diag}\{a_{1},\ldots ,a_{l}\}$. The set of roots are
i) the long ones $\pm \left( \lambda _{i}-\lambda _{j}\right) $ and $\pm
\left( \lambda _{i}+\lambda _{j}\right) $, $1\leq i<j\,\leq l$ and ii) the
short ones $\pm \lambda _{i}$, $1\leq i\leq l$. The set of simple roots is $%
\Sigma =\{\lambda _{1}-\lambda _{2},\ldots ,\lambda _{l-1}-\lambda
_{l},\lambda _{l}\}$, which we write also as $\Sigma =\{\alpha _{1},\ldots
,\alpha _{l}\}$, that is, $\alpha _{i}=\lambda _{i}-\lambda _{i+1}$ if $i<l$
and $\alpha _{l}=\lambda _{l}$.

The Weyl chamber $\mathfrak{a}^{+}\subset \mathfrak{a}$ is defined by the
inequalities 
\begin{equation*}
a_{1}>a_{2}>\cdots >a_{l-1}>a_{l}>0,
\end{equation*}%
and a partial chamber $\mathfrak{a}_{\Theta }\cap \mathrm{cl}\mathfrak{a}%
^{+} $ is defined by a similar relations where some of the strict
inequalities are changed by equalities (e.g. if $\lambda _{i}-\lambda
_{j}\in \Theta $ then $a_{i}=a_{j}$). In particular a characteristic element 
$H_{\Theta }$ for the subset $\Theta =\{\alpha \in \Sigma :\alpha \left(
H_{\Theta }\right) =0\}\subset \Sigma $ is defined by one of these relations.

The subalgebra $\mathfrak{k}$ is composed of the skew-symmetric matrices in $%
\mathfrak{sl}\left( l,l\right) $, that is,%
\begin{equation*}
\left( 
\begin{array}{ccc}
0 & a & a \\ 
-a^{T} & A & B \\ 
-a^{T} & B & A%
\end{array}%
\right) \qquad A+A^{T}=B+B^{T}=0.
\end{equation*}%
It is isomorphic to $\mathfrak{so}\left( l+1\right) \oplus \mathfrak{so}%
\left( l\right) $. The isomorphism is provided by the decomposition%
\begin{equation*}
\left( 
\begin{array}{lll}
0 & a & a \\ 
-a^{T} & A & B \\ 
-a^{T} & B & A%
\end{array}%
\right) =\left( 
\begin{array}{lll}
0 & a & a \\ 
-a^{T} & \left( A+B\right) /2 & \left( A+B\right) /2 \\ 
-a^{T} & \left( A+B\right) /2 & \left( A+B\right) /2%
\end{array}%
\right) +\left( 
\begin{array}{lll}
0 & 0 & 0 \\ 
0 & \left( A-B\right) /2 & -\left( A-B\right) /2 \\ 
0 & -\left( A-B\right) /2 & \left( A-B\right) /2%
\end{array}%
\right) ,
\end{equation*}%
so that $\mathfrak{k}=\mathfrak{k}_{l+1}\oplus \mathfrak{k}_{l}\approx 
\mathfrak{so}\left( l+1\right) \oplus \mathfrak{so}\left( l\right) $ where
the ideals are given by matrices as follows%
\begin{equation}
\mathfrak{k}_{l+1}:\left( 
\begin{array}{ccc}
0 & a & a \\ 
-a^{T} & A & A \\ 
-a^{T} & A & A%
\end{array}%
\right) \qquad \mathfrak{k}_{l}:\left( 
\begin{array}{ccc}
0 & 0 & 0 \\ 
0 & A & -A \\ 
0 & -A & A%
\end{array}%
\right) .  \label{forideaiskbele}
\end{equation}%
In both cases $A$ is skew-symmetric. We write $K_{l+1}=\langle \exp 
\mathfrak{k}_{l+1}\rangle $ and $K_{l+1}=\langle \exp \mathfrak{k}%
_{l}\rangle $.

We start our analysis by describing the irreducible components $V_{\Theta
}^{\sigma }$ defined by the set of roots $\Pi _{\Theta }^{\sigma }$. For
this we separate the cases where $\lambda _{l}$ belongs or not to $\Theta $.

\begin{lema}
\label{lembelecompseml}Suppose that $\lambda _{l}\notin \Theta $ and let 
\begin{equation*}
V_{\Theta }^{\sigma }=\sum_{\alpha \in \Pi _{\Theta }^{\sigma }}\mathfrak{g}%
_{\alpha }
\end{equation*}%
be an irreducible component. Then $\Pi _{\Theta }^{\sigma }$ contains only
short roots or long roots. These sets are described as follows:

\begin{enumerate}
\item \textbf{Short roots:} Take a simple root $\alpha _{i}\notin \Theta $.
Then there are two possibilities:

\begin{enumerate}
\item $\alpha _{i-1}\notin \Theta $. Then $\mathfrak{g}_{-\lambda _{i}}$ is $%
\mathfrak{z}_{\Theta }$-invariant and hence is an irreducible component.

\item $\alpha _{i-1}\in \Theta $. Let $j\left( i\right) <i$ be the smallest
index such that $\{\alpha _{j\left( i\right) },\ldots ,\alpha
_{i-1}\}\subset \Theta $. Then $\Pi _{\Theta }^{\sigma }=\{-\lambda
_{j\left( i\right) },\ldots ,-\lambda _{i-1}\}$ defines a $\mathfrak{z}%
_{\Theta }$-irreducible component.
\end{enumerate}

These sets form a disjoint union of the negative short roots $-\lambda _{i}$%
, $1\leq i\leq l$. (Note that this disjoint union completely determines $%
\Theta $.)

\item \textbf{Long roots:} A subset $\Pi _{\Theta }^{\sigma }$ contains only
roots of the type $\lambda _{i}-\lambda _{j}$ or of the type $-\lambda
_{i}-\lambda _{j}$.
\end{enumerate}
\end{lema}

\begin{profe}
To see the components corresponding to the short roots take an index $i$
with $\alpha _{i}\notin \Theta $. An easy check shows that the only simple
roots $\alpha $ such that $-\lambda _{i}+\alpha $ is a root are $\alpha
=\lambda _{i}-\lambda _{i+1}$ or $\alpha =\lambda _{l}$. By assumption these
simple roots are not in $\Theta $. This implies that $-\lambda _{i}$ is the
highest weight of an irreducible representation of $\mathfrak{g}\left(
\Theta \right) $. The weights of this representation are restrictions of $%
\mathfrak{a}\left( \Theta \right) $ of roots. They have the form $-\lambda
_{i}-\beta _{1}-\cdots -\beta _{k}$ with $\beta _{i}\in \Theta $. But these
successive diferences are roots only when $\beta _{1}=\lambda _{i-1}-\lambda
_{i}$, $\beta _{2}=\lambda _{i-2}-\lambda _{i-1}$, and so on, obtaining the $%
-\lambda _{i}$, $-\lambda _{i-1}$, \ldots , $-\lambda _{j\left( i\right) }$
with $j\left( i\right) $ as in (b). This concludes the case of the short
roots.

Now, take a long root, e.g. $\lambda _{i}-\lambda _{j}$. Then $\Pi _{\Theta
}\left( \lambda _{i}-\lambda _{j}\right) $ does not contain short roots that
were already exhausted. On the other hand, by assumption $\Theta $ is
contained in the set of roots of the type $\lambda _{r}-\lambda _{s}$. Since
this set is closed by sum we conclude that the roots in $\Pi _{\Theta
}\left( \lambda _{i}-\lambda _{j}\right) $ have the type $\lambda
_{r}-\lambda _{s}$. The same argument applies to $\Pi _{\Theta }\left(
-\lambda _{i}-\lambda _{j}\right) $.
\end{profe}

\begin{lema}
\label{lembelecompcoml}Suppose that $\lambda _{l}\in \Theta $ and let $i_{0}$
be the largest index such that $\lambda _{i_{0}}-\lambda _{i_{0}+1}\notin
\Theta $, that is, $\{\alpha _{i_{0}+1},\ldots ,\alpha _{l}=\lambda _{l}\}$
is the connected component of $\Theta $ containing $\lambda _{l}$. Then the
sets of roots defining the $\mathfrak{z}_{\Theta }$-irreducible components
are as follows:

\begin{enumerate}
\item \textbf{Components containing short roots}: If $i\leq i_{0}$ then $\Pi
_{\Theta }\left( -\lambda _{i}\right) $ contains $-\lambda _{i}+\lambda _{k}$
and $-\lambda _{i}-\lambda _{k}$ for all $k\geq i_{0}+1$. (The short roots $%
-\lambda _{i}$, $i>i_{0}$, belong to $\langle \Theta \rangle ^{-}$.)

Moreover the sets of short roots belonging to the same component are as in
Lemma \ref{lembelecompseml} (1), namely $\{-\lambda _{j\left( i\right)
},\ldots ,-\lambda _{i-1}\}$ where $\{\alpha _{j\left( i\right) },\ldots
,\alpha _{i-1}\}$ is a connected component of $\Theta $.

\item \textbf{Components containing only long roots}: If $i<j\leq i_{0}$
then $\Pi _{\Theta }\left( -\lambda _{i}+\lambda _{j}\right) $ has only
roots $\lambda _{r}-\lambda _{s}$ and $\Pi _{\Theta }\left( -\lambda
_{i}-\lambda _{j}\right) $ has only roots $-\lambda _{r}-\lambda _{s}$.
\end{enumerate}

(These sets exhaust the roots because $-\lambda _{i},-\lambda _{i}\pm
\lambda _{j}\in \langle \Theta \rangle ^{-}$ if $i_{0}+1\leq i<j$.)
\end{lema}

\begin{profe}
By assumption $\Theta $ contains the subdiagram simple roots $%
B_{l-i_{0}}=\{\alpha _{i_{0}+1},\ldots ,\alpha _{l}\}$. This implies that
the roots $\pm \lambda _{k}\pm \lambda _{j}$ and $\pm \lambda _{k}$ belong
to $\langle \Theta \rangle $ if $i_{0}+1\leq k<j$. Take a short root $%
-\lambda _{i}$ with $i\leq i_{0}$, which is not in $\langle \Theta \rangle
^{-}$. For any root $\alpha \in \langle \Theta \rangle $ such that $-\lambda
_{i}+\alpha $ is a root we have $-\lambda _{i}+\alpha \in \Pi _{\Theta
}\left( -\lambda _{i}\right) $. If we take $\alpha =\pm \lambda _{k}$, $%
k\geq i_{0}+1$, we see that $-\lambda _{i}\pm \lambda _{k}\in \Pi _{\Theta
}\left( -\lambda _{i}\right) $, proving the first part of (1). By the same
argument of the proof Lemma \ref{lembelecompseml} we get the statement about
the short roots.

Now, a long root $-\lambda _{i}+\lambda _{j}$, $i<j\leq i_{0}$, is
orthogonal to every root in $B_{l-i_{0}}$. Hence the only way to get new
roots from $-\lambda _{i}+\lambda _{j}$ is by adding or subtracting roots in 
$\Theta \setminus B_{l-i_{0}}$. These roots have the type $\lambda
_{r}-\lambda _{s}$, so that as in the proof of Lemma \ref{lembelecompseml}
we see that $\Pi _{\Theta }\left( -\lambda _{i}+\lambda _{j}\right) $
contains only roots of the type $\lambda _{r}-\lambda _{s}$. The same
argument works for $-\lambda _{i}-\lambda _{j}$, $i<j\leq i_{0}$, showing
(2).
\end{profe}

The next step is to look at the $K_{\Theta }$-irreducibility of the $%
\mathfrak{z}_{\Theta }$-irreducible components. For this we use the $M$%
-equivalence classes so we are led to consider separetely different values
of $l$.

\begin{lema}
Take $B_{l}$ with $l\geq 5$.

\begin{enumerate}
\item Suppose that $\lambda _{l}\notin \Theta $. Then any component $%
V_{\Theta }^{\sigma }$ is $K_{\Theta }$-irreducible.

\item If $\lambda _{l}\in \Theta $ then $K_{\Theta }$ is irreducible in the
components $V_{\Theta }^{\sigma }$ such that $\Pi _{\Theta }^{\sigma }$
contains only long roots as in Lemma \ref{lembelecompcoml} (2).
\end{enumerate}
\end{lema}

\begin{profe}
We just piece together different facts proved previously. First if $\lambda
_{l}\notin \Theta $ then $\Theta $ contains only long roots. Hence by Lemma %
\ref{lemtransimplelong} the subgroup $\mathcal{W}_{\Theta }$ acts
transitively on the sets $\Pi _{\Theta }^{\sigma }$ for any irreducible
component $V_{\Theta }^{\sigma }$. Now if $l\geq 5$ then the $M$-equivalence
classes are the short roots $\{-\lambda _{i}\}$ and $\{-\lambda _{i}+\lambda
_{j},-\lambda _{i}-\lambda _{j}\}$, $i<j$. Hence by Lemma \ref%
{lembelecompseml} the intersection of a $M$-equivalence class with a set $%
\Pi _{\Theta }^{\sigma }$ has just one root. Therefore, the assumptions of
Lemma \ref{lemirreducible} are satisfied, and we get the conclusion that any 
$V_{\Theta }^{\sigma }$ is $K_{\Theta }$-irreducible, proving (1).

The proof of (2) is similar. Take a subset $\Pi _{\Theta }^{\sigma }$ as in
the statement. Again no two roots in $\Pi _{\Theta }^{\sigma }$ are $M$%
-equivalent. As to the transitive action of $\mathcal{W}_{\Theta }$ consider
the subset $B_{l-i_{0}}=\{\alpha _{i_{0}+1},\ldots ,\alpha _{l}\}$ defined
in the proof of Lemma \ref{lembelecompcoml}. Then $\mathcal{W}_{\Theta }$ is
the direct product $\mathcal{W}_{\Theta }=\mathcal{W}_{\Theta \setminus
B_{l-i_{0}}}\times \mathcal{W}_{B_{l-i_{0}}}$, and any $w\in \mathcal{W}%
_{B_{l-i_{0}}}$ is the identity in $\Pi _{\Theta }^{\sigma }$, because the
sets $\Pi _{\Theta }^{\sigma }$ and $B_{l-i_{0}}$ are orthogonal (see the
proof of Lemma \ref{lembelecompcoml}). Now $\mathcal{W}_{\Theta \setminus
B_{l-i_{0}}}$ acts transitively on $\Pi _{\Theta }^{\sigma }$ since $\Theta
\setminus B_{l-i_{0}}$ has only long roots (see the proof of Lemma \ref%
{lemtransimplelong}).
\end{profe}

It remains to analyze the components $V_{\Theta }\left( -\lambda _{i}\right) 
$ containing short roots $-\lambda _{i}$ in case $\lambda _{l}\in \Theta $.
Contrary to the others these are not $K_{\Theta }$-irreducible. Let us write
them explicitly as follows: Let $i_{0}$ be, as in Lemma \ref{lembelecompcoml}%
, the largest index such that $\alpha _{i_{0}}=\lambda _{i_{0}}-\lambda
_{i_{0}+1}\notin \Theta $, so that $-\lambda _{i}\notin \langle \Theta
\rangle ^{-}$. For $i\leq i_{0}$ and $k>i_{0}$ put 
\begin{equation*}
W_{\Theta }^{ik}=\mathfrak{g}_{-\lambda _{i}+\lambda _{k}}\oplus \mathfrak{g}%
_{-\lambda _{i}}\oplus \mathfrak{g}_{-\lambda _{i}-\lambda _{k}}\qquad 
\mathrm{and}\qquad W_{\Theta }^{i}=\sum_{k\geq i_{0}+1}W_{\Theta }^{ik}.
\end{equation*}%
(Note that the last sum is not direct.)

By Lemma \ref{lembelecompcoml} (2) we have one irreducible component $%
V_{\Theta }\left( -\lambda _{i}\right) $ for each index $i\leq i_{0}$ such
that $\alpha _{i}=\lambda _{i}-\lambda _{i+1}\notin \Theta $. To write it in
terms of the subspaces $W_{\Theta }^{i}$ let $j\left( i\right) \leq i$ be
defined by

\begin{enumerate}
\item $j\left( i\right) =i$ if $\alpha _{i-1}\notin \Theta $, and

\item $j\left( i\right) $ is such that $\{\alpha _{j\left( i\right) },\ldots
,\alpha _{i-1}\}$ is a connected component of $\Theta $ if $\alpha _{i-1}\in
\Theta $.
\end{enumerate}

Then 
\begin{equation*}
V_{\Theta }\left( -\lambda _{i}\right) =\sum_{k=j\left( i\right)
}^{i}W_{\Theta }^{k}.
\end{equation*}

Now for $i,j$ let 
\begin{equation}
E_{ij}^{-}=\left( 
\begin{array}{ccc}
0 & 0 & 0 \\ 
0 & E_{ij} & 0 \\ 
0 & 0 & -E_{ji}%
\end{array}%
\right) ~E_{ij}^{+}=\left( 
\begin{array}{ccc}
0 & 0 & 0 \\ 
0 & 0 & 0 \\ 
0 & E_{ij}-E_{ji} & 0%
\end{array}%
\right) ~E_{i}^{0}=\left( 
\begin{array}{ccc}
0 & e_{i} & 0 \\ 
0 & 0 & 0 \\ 
-e_{i}^{T} & 0 & 0%
\end{array}%
\right)  \label{forbasesbele}
\end{equation}%
where $E_{ij}$ and $e_{i}$ are basic $l\times l$ and $1\times l$ matrices.
These matrices are generators of $\mathfrak{g}_{-\lambda _{i}+\lambda _{j}}$%
, $\mathfrak{g}_{-\lambda _{i}-\lambda _{j}}$ and $\mathfrak{g}_{-\lambda
_{i}}$, respectively. So that $\{E_{ik}^{-},E_{ik}^{+},E_{ik}^{0}\}$ is a
basis of $W_{\Theta }^{ik}$ and $\{E_{ik}^{-},E_{ik}^{+},E_{ik}^{0}:k\geq
i_{0+1}\}$ is a basis of $W_{\Theta }^{i}$.

Before proceeding we note that the subspace $W_{\Theta }^{ik}$ is invariant
and irreducible by adjoint representation of the subalgebra $\mathfrak{g}%
\left( \lambda _{k}\right) \approx \mathfrak{sl}\left( 2,\mathbb{R}\right) $
generated by $\mathfrak{g}_{\pm \lambda _{k}}$, thus defining an irreducible
representation of $\mathfrak{sl}\left( 2,\mathbb{R}\right) $. By
dimensionality this representation is equivalent to the adjoint
representation of $\mathfrak{sl}\left( 2,\mathbb{R}\right) $, which in turn
is not $\mathfrak{so}\left( 2\right) $-irreducible: It decomposes into the
subspaces of skew-symmetric ($1$-dimensional) and symmetric ($2$%
-dimensional) matrices. The equivalence $\mathfrak{g}\left( \lambda
_{k}\right) \approx \mathfrak{sl}\left( 2,\mathbb{R}\right) $ maps $%
\mathfrak{g}_{-\lambda _{i}+\lambda _{k}}$ and $\mathfrak{g}_{-\lambda
_{i}-\lambda _{k}}$ onto the upper and lower triangular matrices,
respectively and $\mathfrak{g}_{-\lambda _{i}}$ onto the diagonal matrices.
So we get

\begin{lema}
The representation of $\mathfrak{g}\left( \lambda _{k}\right) $ decomposes $%
W_{\Theta }^{ik}$ into two $\mathfrak{k}_{\{\lambda _{k}\}}$-invariant
subspaces, namely%
\begin{equation*}
\left( W_{\Theta }^{ik}\right) _{1}=\mathrm{span}\{E_{ik}^{-}-E_{ik}^{+}\}%
\qquad \mathrm{and}\qquad \left( W_{\Theta }^{ik}\right) _{2}=\mathrm{span}%
\{E_{ik}^{-}+E_{ik}^{+},E_{ik}^{0}\}.
\end{equation*}
\end{lema}

Now we can decompose $V_{\Theta }\left( -\lambda _{i}\right) $, $i\leq i_{0}$
when $\lambda _{l}\in \Theta $ into $K_{\Theta }$-irreducible subspaces.

\begin{lema}
\label{lemirreducw}In $B_{l}$, $l\geq 5$, suppose $\lambda _{l}\in \Theta $
and let $i_{0}$ be the largest index such that $\lambda _{i_{0}}-\lambda
_{i_{0}+1}\notin \Theta $. If $i\leq i_{0}$ then $-\lambda _{i}\notin
\langle \Theta \rangle ^{-}$ and the $\mathfrak{z}_{\Theta }$-irreducible
component $V_{\Theta }\left( -\lambda _{i}\right) $ is the direct sum of the
following $K_{\Theta }$-irreducible and invariant subspaces 
\begin{equation*}
\left( W_{\Theta }^{i}\right) _{1}=\sum_{j=j\left( i\right) }^{i}\sum_{k\geq
i_{0}+1}\left( W_{\Theta }^{jk}\right) _{1}\qquad \mathrm{and}\qquad \left(
W_{\Theta }^{i}\right) _{2}=\sum_{j=j\left( i\right) }^{i}\sum_{k\geq
i_{0}+1}\left( W_{\Theta }^{jk}\right) _{2}
\end{equation*}%
with $\dim \left( W_{\Theta }^{i}\right) _{1}=l-i_{0}+i-j\left( i\right) +1$
and $\dim \left( W_{\Theta }^{i}\right) _{2}=2\left( l-i_{0}+i-j\left(
i\right) +1\right) $.

Furthermore, $\left( W_{\Theta }^{i}\right) _{1}=V_{\Theta }\left( -\lambda
_{i}\right) \cap T_{b_{\Theta }}K_{l}\cdot b_{\Theta }$ and $\left(
W_{\Theta }^{i}\right) _{2}=V_{\Theta }\left( -\lambda _{i}\right) \cap
T_{b_{\Theta }}K_{l+1}\cdot b_{\Theta }$.
\end{lema}

\begin{profe}
The intersections in the last statement with the tangent space to the orbits 
$K_{l}\cdot b_{\Theta }$ and $K_{l+1}\cdot b_{\Theta }$ are readily obtained
from the matrices in $\mathfrak{k}_{l}$ and $\mathfrak{k}_{l+1}$ given in (%
\ref{forideaiskbele}) and the definition of the subspaces $\left( W_{\Theta
}^{i}\right) _{1}$ and $\left( W_{\Theta }^{i}\right) _{2}$. It follows by
Lemma \ref{lemorbitnormal} that these subspaces are $K_{\Theta }$-invariant.

To check irreducibility consider $\left( W_{\Theta }^{i}\right) _{2}$ and
take a nonzero $K_{\Theta }$-invariant subspace $Z\subset \left( W_{\Theta
}^{i}\right) _{2}$.

We claim that there are $j\in \left[ j\left( i\right) ,i\right] $ and $k\geq
i_{0}+1$ such that $\left( W_{\Theta }^{jk}\right) _{2}\subset Z$.

By Lemma \ref{leminterinvMclass} we have a nontrivial intersection of $Z$
with a subspace $\sum_{\alpha }\mathfrak{g}_{\alpha }$, with the sum
extended to a $M$-equivalence class. Since we are assuming that $l\geq 5$,
the $M$-equivalence classes are $\{-\lambda _{s}+\lambda _{r},-\lambda
_{s}-\lambda _{r}\}$ and $\{-\lambda _{s}\}$. Hence either there exists $%
j\in \left[ j\left( i\right) ,i\right] $ such that $Z\cap \mathfrak{g}%
_{-\lambda _{j}}\neq \{0\}$ or there are $j\in \left[ j\left( i\right) ,i%
\right] $ and $k\geq i_{0}+1$ such that $Z\cap \left( \mathfrak{g}_{-\lambda
_{j}+\lambda _{k}}\oplus \mathfrak{g}_{-\lambda _{j}-\lambda _{k}}\right)
\neq \{0\}$. In both cases we have 
\begin{equation*}
Z\cap \left( W_{\Theta }^{jk}\right) _{2}\neq \{0\}
\end{equation*}%
However, $\left( W_{\Theta }^{jk}\right) _{2}$ is invariant and irreducible
for $\mathfrak{k}_{\{\lambda _{k}\}}$. Since $k\geq i_{0}+1$ we have $%
\lambda _{k}\in \langle \Theta \rangle $ and $\mathfrak{k}_{\{\lambda
_{k}\}}\subset \mathfrak{k}_{\Theta }$. By $K_{\Theta }$-invariance of $Z$
we conclude that $\left( W_{\Theta }^{jk}\right) _{2}\subset Z$.

Now let $B_{l-i_{0}+1}$ be the connected component of $\Theta $ containing $%
\lambda _{l}$. Then its Weyl group $\mathcal{W}_{B_{l-i_{0}+1}}\subset 
\mathcal{W}_{\Theta }$ acts transitively on the set of its short root. This
means that if $k_{1},k_{2}\geq i_{0}+1$ then there exists $w\in \mathcal{W}%
_{B_{l-i_{0}+1}}$ such that $w\lambda _{k_{1}}=\lambda _{k_{2}}$. Combining
this transitivity with the claim it follows that $\left( W_{\Theta
}^{jk}\right) _{2}\subset Z$ for every $k\geq i_{0}+1$. Consequently, there
exists $j\in \left[ j\left( i\right) ,i\right] $ such that $\sum_{k\geq
i_{0}+1}\left( W_{\Theta }^{jk}\right) _{2}\subset Z$.

To finish the proof we use the subgroup $\mathcal{W}_{\left[ j\left(
i\right) ,i-1\right] }$ of $\mathcal{W}_{\Theta }$ generated by the
reflections with respect to the roots in the connected component $\{\alpha
_{j\left( i\right) },\ldots ,\alpha _{i-1}\}\subset \Theta $. This subgroup
is the permutation group of $\{j\left( i\right) ,\ldots ,i-1,i\}$. Since $%
\sum_{k\geq i_{0}+1}\left( W_{\Theta }^{jk}\right) _{2}\subset Z$ for some $%
j\in \left[ j\left( i\right) ,i\right] $ we conclude $\sum_{k\geq
i_{0}+1}\left( W_{\Theta }^{sk}\right) _{2}\subset Z$ for all $s\in \left[
j\left( i\right) ,i\right] $, so that $\left( W_{\Theta }^{i}\right)
_{2}\subset Z$, showing irreducibility of $\left( W_{\Theta }^{i}\right)
_{2} $.

The proof for $\left( W_{\Theta }^{i}\right) _{1}$ is similar.
\end{profe}

Summarizing the above discussion we have the following $K_{\Theta }$%
-invariant irreducible subspaces for $B_{l}$, $l\geq 5$:

\begin{enumerate}
\item A $\mathfrak{z}_{\Theta }$-component $V_{\Theta }\left( -\lambda
_{i}\right) $ containing only short roots. These components occur only when $%
\lambda _{l}\notin \Theta $.

\item $\mathfrak{z}_{\Theta }$-components \ $V_{\Theta }\left( -\lambda
_{i}+\lambda _{j}\right) $ and $V_{\Theta }\left( -\lambda _{i}-\lambda
_{j}\right) $, $i<j$, containing only long roots. These subspaces occur in
both cases when $\lambda _{l}$ belongs or not to $\Theta $. When $\lambda
_{l}\in \Theta $ the indexes $i,j$ satisfy $i<j\leq i_{0}$ where $\{\alpha
_{i_{0}},\ldots ,\alpha _{l}=\lambda _{l}\}$ is the connected component of $%
\Theta $ containing $\lambda _{l}$.

\item The subspaces $\left( W_{\Theta }^{i}\right) _{1}$ and $\left(
W_{\Theta }^{i}\right) _{2}$ contained in a $\mathfrak{z}_{\Theta }$%
-component \ $V_{\Theta }\left( -\lambda _{i}\right) $ when $\lambda _{l}\in
\Theta $.
\end{enumerate}

These are not the only invariant irreducible subspaces of $K_{\Theta }$,
since among them some pairs $V_{1}\neq V_{2}$ are $K_{\Theta }$-equivalent,
enabling the existence of invariant subspaces inside $V_{1}\oplus V_{2}$.
Among these pairs we can discard the following by $M$-equivalence we discard
the following pairs: i) $V_{1}$ is a subspace in item (1) and $V_{2}$ in (1)
or (2); ii) $V_{1}$ is a subspace in (2) and $V_{2}$ in (3); iii) $V_{1}$ is
a subspace $\left( W_{\Theta }^{i}\right) _{1,2}$ and $V_{2}=\left(
W_{\Theta }^{j}\right) _{1,2}$ if $i\neq j$. Since (1) and (3) are subspaces
for different $\Theta $ and $\left( W_{\Theta }^{i}\right) _{1}$ and $\left(
W_{\Theta }^{i}\right) _{2}$ do not have the same dimension, it remains the
subspaces $V_{\Theta }\left( -\lambda _{i}+\lambda _{j}\right) $ and $%
V_{\Theta }\left( -\lambda _{i}-\lambda _{j}\right) $ of (2). These are
indeed equivalent.

\begin{lema}
In $B_{l}$, $l\geq 5$, the subspaces $V_{\Theta }\left( -\lambda
_{i}+\lambda _{j}\right) $ and $V_{\Theta }\left( -\lambda _{i}-\lambda
_{j}\right) $ as in (2) above are $K_{\Theta }$-equivalent if both roots $%
-\lambda _{i}+\lambda _{j}$ and $-\lambda _{i}-\lambda _{j}$ do not belong
to $\langle \Theta \rangle ^{-}$.
\end{lema}

\begin{profe}
To prove equivalence we shall exhibit a proper $K_{\Theta }$-invariant
subspace $\{0\}\neq V\subset V_{\Theta }\left( -\lambda _{i}+\lambda
_{j}\right) \oplus V_{\Theta }\left( -\lambda _{i}-\lambda _{j}\right) $
different from the irreducible components $V_{\Theta }\left( -\lambda
_{i}+\lambda _{j}\right) $ and $V_{\Theta }\left( -\lambda _{i}-\lambda
_{j}\right) $. This implies equivalence by Proposition \ref{propsumirredcomp}%
).

The required subspace $V$ will be obtained from the tangent space at the
origin of the orbit of the normal subgroup $K_{l}$. By (\ref{forideaiskbele}%
) the matrices in the Lie algebra $\mathfrak{k}_{l}$ of $K_{l}$ are 
\begin{equation*}
\left( 
\begin{array}{ccc}
0 & 0 & 0 \\ 
0 & A & -A \\ 
0 & -A & A%
\end{array}%
\right) \quad A+A^{T}=0.
\end{equation*}%
Looking at these matrices we see that after identifying $T_{b_{\Theta }}%
\mathbb{F}_{\Theta }$ with $\mathfrak{n}_{\Theta }^{-}$ the tangent space $%
T_{b_{\Theta }}\left( K_{l}\cdot b_{\Theta }\right) $ is identified to the
subspace $T_{l}\subset \mathfrak{n}_{\Theta }^{-}$ spanned by $\mathrm{pr}%
\left( E_{rs}^{-}-E_{rs}^{+}\right) $, $r>s$, where $E_{rs}^{\pm }$ were
defined in (\ref{forbasesbele}) and $\mathrm{pr}:\mathfrak{n}^{-}\rightarrow 
\mathfrak{n}_{\Theta }^{-}$ is the projection w.r.t. the root spaces
decomposition.

The tangent space $T_{b_{\Theta }}\left( K_{l}\cdot b_{\Theta }\right) $ is
invariant by the isotropy representation of $K_{\Theta }$, by Lemma \ref%
{lemorbitnormal}. Hence $T_{l}$ is invariant by the adjoint action of $%
K_{\Theta }$.

Now if $-\lambda _{r}+\lambda _{s}$, $-\lambda _{r}-\lambda _{s}\in \Pi
^{-}\setminus \langle \Theta \rangle ^{-}$ then $E_{rs}^{-}-E_{rs}^{+}=%
\mathrm{pr}\left( E_{rs}^{-}-E_{rs}^{+}\right) $. Hence the following
vectors form a basis of $T_{l}$:

\begin{enumerate}
\item $E_{rs}^{-}$ such that $-\lambda _{r}+\lambda _{s}\in \Pi _{\Theta
}\left( -\lambda _{i}+\lambda _{j}\right) $ and $-\lambda _{r}-\lambda
_{s}\notin \Pi _{\Theta }\left( -\lambda _{i}-\lambda _{j}\right) $.

\item $E_{rs}^{+}$ such that $-\lambda _{r}-\lambda _{s}\in \Pi _{\Theta
}\left( -\lambda _{i}-\lambda _{j}\right) $ and $-\lambda _{r}+\lambda
_{s}\notin \Pi _{\Theta }\left( -\lambda _{i}+\lambda _{j}\right) $.

\item $E_{rs}^{-}-E_{rs}^{+}$ such that $-\lambda _{r}+\lambda _{s}\in \Pi
_{\Theta }\left( -\lambda _{i}+\lambda _{j}\right) $ and $-\lambda
_{r}-\lambda _{s}\in \Pi _{\Theta }\left( -\lambda _{i}-\lambda _{j}\right) $%
.
\end{enumerate}

The third case is not empty (e.g. $\left( r,s\right) =\left( i,j\right) $
fall in this case), which means that $E_{rs}^{-}-E_{rs}^{+}\in T_{l}$ for
some pair $\left( r,s\right) $. For this pair $T_{l}\cap V_{\Theta }\left(
\lambda _{r}-\lambda _{s}\right) =T_{1}\cap V_{\Theta }\left( -\lambda
_{r}-\lambda _{s}\right) =\{0\}$, which shows that $T_{l}$ is proper and
different from $V_{\Theta }\left( -\lambda _{i}+\lambda _{j}\right) $ and $%
V_{\Theta }\left( -\lambda _{i}-\lambda _{j}\right) $. By Proposition \ref%
{propsumirredcomp} it follows that these irreducible subspaces are $%
K_{\Theta }$-equivalent.
\end{profe}

In conclusion we have:

\begin{teorema}
Let $\mathbb{F}_{\Theta }$ be a flag manifold of $B_{l}=\mathfrak{so}\left(
l+1,l\right) $, $l\geq 5$. Then the $K_{\Theta }$-invariant irreducible
subspaces of $\mathfrak{n}_{\Theta }^{-}$ are in the following classes:

\begin{enumerate}
\item Isolated subspaces:

\begin{enumerate}
\item $V_{\Theta }\left( -\lambda _{i}\right) $ when $\lambda _{l}\notin
\Theta $. These subspaces contain root spaces of short roots only.

\item $V_{\Theta }\left( -\lambda _{i}+\lambda _{j}\right) $, $i<j$, when $%
-\lambda _{i}-\lambda _{j}\in \langle \Theta \rangle ^{-}$. Any such pair
occur if $-\lambda _{l}\in \Theta $. Otherwise we have $i<j\leq i_{0}$,
where $\{\alpha _{i_{0}+1},\ldots ,\alpha _{l}=\lambda _{l}\}$ is the
connected component of $\Theta $ containing $\lambda _{l}$.

\item The same as (b) interchanging the roles of $-\lambda _{i}+\lambda _{j}$
and $-\lambda _{i}-\lambda _{j}$.

\item The subspaces%
\begin{equation*}
\left( W_{\Theta }^{i}\right) _{1}=\sum_{j=j\left( i\right) }^{i}\sum_{k\geq
i_{0}+1}\left( W_{\Theta }^{jk}\right) _{1}\qquad \mathrm{and}\qquad \left(
W_{\Theta }^{i}\right) _{2}=\sum_{j=j\left( i\right) }^{i}\sum_{k\geq
i_{0}+1}\left( W_{\Theta }^{jk}\right) _{2}
\end{equation*}%
defined in Lemma \ref{lemirreducw}. These subspaces decompose $V\left(
-\lambda _{i}\right) $ when $\lambda _{l}\in \Theta $.
\end{enumerate}

\item A continuum of invariant subspaces parametrized by $[(x,y)]\in \mathbb{%
R}P^{2}$ given by 
\begin{equation*}
V_{[(x,y)]}^{ij}=\left\{ xX+yTX:X\in V_{\Theta }\left( -\lambda _{i}+\lambda
_{j}\right) \right\} \qquad i<j.
\end{equation*}%
The indexes $ij$ are as in (1.b), and here both $-\lambda _{i}+\lambda _{j}$
and $-\lambda _{i}-\lambda _{j}$ are not in $\langle \Theta \rangle ^{-}$.
\end{enumerate}
\end{teorema}

The low dimensional cases $l=2,3,4$ must be treated separetely because of
the difference in the $M$-equivalence classes.

For instance $B_{2}$ has three flag manifolds. The maximal one whose
irreducible components are detected by the $M$-equivalence classes $%
\{\lambda _{1}-\lambda _{2},\lambda _{1}+\lambda _{2}\}$ and $\{\lambda
_{1},\lambda _{2}\}$. Hence, there are two continuous families of $1$%
-dimensional irreducible subspaces. In the flag $\mathbb{F}_{\{\lambda
_{1}-\lambda _{2}\}}$ there are two $\mathfrak{z}_{\Theta }$-irreducible
components defined by the sets $\{-\lambda _{2},-\lambda _{1}\}$ and $%
\{-\lambda _{1}-\lambda _{2}\}$. Both are $K_{\Theta }$-irreducible and
clearly they are not equivalent. On the other hand the flag $\mathbb{F}%
_{\{\lambda _{2}\}}$ has just one $3$-dimensional $\mathfrak{z}_{\Theta }$%
-irreducible component which decomposes into a $1$-dimensional plus a $2$%
-dimensional irreducible subspaces of $K_{\Theta }$ (as happens to the
adjoint representation of $\mathfrak{sl}\left( 2,\mathbb{R}\right) $).

For $B_{3}$ and $B_{4}$ the compact subalgebra $\mathfrak{k}$ ($\mathfrak{so}%
\left( 3\right) \oplus \mathfrak{so}\left( 4\right) $ and $\mathfrak{so}%
\left( 4\right) \oplus \mathfrak{so}\left( 5\right) $, respectively) splits
once more because $\mathfrak{so}\left( 4\right) =\mathfrak{so}\left(
3\right) \oplus \mathfrak{so}\left( 3\right) $. By Lemma \ref{lemorbitnormal}
these simple components of $\mathfrak{k}$ can yield new $K_{\Theta }$%
-invariant subspaces. invariant subspaces. The example with $D_{4}$ below,
which has a similar splitting, illustrates this occurence of new invariant
subspaces. Another aspect that differs $B_{3}$ and $B_{4}$ from the general
case are the $M$-equivalence classes that have more elements. This can
introduce more $K_{\Theta }$-equivalence than the general case. The example
with $C_{4}$ below illustrates this fact.

\subsection{Flags of $C_{l}=\mathfrak{sp}\left( l,\mathbb{R}\right) $}

The symplectic Lie algebra $\mathfrak{sp}\left( l,\mathbb{R}\right) $ is
composed of the real $2l\times 2l$ matrices%
\begin{equation*}
\left( 
\begin{array}{cc}
A & B \\ 
C & -A^{T}%
\end{array}%
\right) \qquad B-B^{T}=C-C^{T}=0
\end{equation*}%
written in the basis $\{e_{1},\ldots ,e_{l},f_{1},\ldots ,f_{l}\}$. In this
case $\mathfrak{a}$ is the subalgebra of matrices 
\begin{equation*}
\left( 
\begin{array}{ll}
\Lambda & 0 \\ 
0 & -\Lambda%
\end{array}%
\right)
\end{equation*}%
with $\Lambda =\mathrm{diag}\{a_{1},\ldots ,a_{l}\}$. The set of roots are
i) the long ones $\pm 2\lambda _{i}$, $1\leq i\leq l$ and ii) the short ones 
$\pm \left( \lambda _{i}-\lambda _{j}\right) $ and $\pm \left( \lambda
_{i}+\lambda _{j}\right) $, $1\leq i<j\,\leq l$. The set of simple roots is $%
\Sigma =\{\lambda _{1}-\lambda _{2},\ldots ,\lambda _{l-1}-\lambda
_{l},2\lambda _{l}\}$, which we write also as $\Sigma =\{\alpha _{1},\ldots
,\alpha _{l}\}$, that is, $\alpha _{i}=\lambda _{i}-\lambda _{i+1}$ if $i<l$
and $\alpha _{l}=2\lambda _{l}$.

The Weyl chamber $\mathfrak{a}^{+}\subset \mathfrak{a}$ is defined by the
inequalities 
\begin{equation*}
a_{1}>a_{2}>\cdots >a_{l-1}>a_{l}>0,
\end{equation*}%
and a partial chamber $\mathfrak{a}_{\Theta }\cap \mathrm{cl}\mathfrak{a}%
^{+} $ is defined by a similar relations where some of the strict
inequalities are changed by equalities (e.g. if $\lambda _{i}-\lambda
_{j}\in \Theta $ then $a_{i}=a_{j}$). In particular a characteristic element 
$H_{\Theta }$ for the subset $\Theta =\{\alpha \in \Sigma :\alpha \left(
H_{\Theta }\right) =0\}\subset \Sigma $ is defined by one of these relations.

The subalgebra $\mathfrak{k}$ is composed of the skew-symmetric matrices in $%
\mathfrak{sp}\left( l,\mathbb{R}\right) $, that is,%
\begin{equation*}
\left( 
\begin{array}{cc}
A & -B \\ 
B & A%
\end{array}%
\right) \qquad A+A^{T}=B-B^{T}=0.
\end{equation*}%
It is isomorphic to $\mathfrak{u}\left( l\right) =\mathfrak{su}\left(
l\right) \oplus \mathbb{R}$, where the isomorphism associates the above
matrix the complex matrix $A+iB$.

To describe the $\mathfrak{z}_{\Theta }$-irreducible components $V_{\Theta
}^{\sigma }$ defined by the set of roots $\Pi _{\Theta }^{\sigma }$ we
consider first the components $V_{\Theta }\left( -2\lambda _{i}\right) $
containing the long roots.

\begin{lema}
\label{lemcelelong}Let $i=1,\ldots ,l$ be an index such that $\alpha
_{i}\notin \Theta $.

\begin{enumerate}
\item If $i=1$ or $\alpha _{i-1}\notin \Theta $ then $V_{\Theta }\left(
-2\lambda _{i}\right) =\mathfrak{g}_{-2\lambda _{i}}$.

\item Otherwise let $j\left( i\right) <i$ be such that $\{\alpha _{j\left(
i\right) },\ldots ,\alpha _{i-1}\}$ is the connected component of $\Theta $
containing $\alpha _{i-1}$. Then 
\begin{equation}
V_{\Theta }\left( -2\lambda _{i}\right) =\sum_{k,r=j\left( i\right) }^{i}%
\mathfrak{g}_{-\lambda _{k}-\lambda _{r}}.  \label{forcelelong}
\end{equation}
\end{enumerate}

If $\alpha _{i}\in \Theta $ then either $-2\lambda _{i}\in \langle \Theta
\rangle $ if $2\lambda _{l}\in \Theta $ and $\alpha _{i}$ and $2\lambda _{l} 
$ are in the same connected component of $\Theta $ or $-2\lambda _{i}\in \Pi
_{\Theta }\left( -2\lambda _{j}\right) $ where $j>i$ is the smallest index
such that $\alpha _{j}\notin \Theta $.
\end{lema}

\begin{profe}
Since $\alpha _{i}=\lambda _{i}-\lambda _{i+1}\notin \Theta $ the only way
that $-2\lambda _{i}\pm \alpha $ is a root with $\alpha \in \Theta $ is in
the string $-2\lambda _{i}-\left( \lambda _{i-1}-\lambda _{i}\right)
=-\lambda _{i-1}-\lambda _{i}$ and $-2\lambda _{i}-2\left( \lambda
_{i-1}-\lambda _{i}\right) =-2\lambda _{i-1}$. Hence if $\alpha _{i-1}\notin
\Theta $ (or $i=1$) no such sum occurs and $V_{\Theta }\left( -2\lambda
_{i}\right) =\mathfrak{g}_{-2\lambda _{i}}$.

On the other hand if $\alpha _{i-1}=\lambda _{i-1}-\lambda _{i}\in \Theta $
then the roots $-\lambda _{i-1}-\lambda _{i}=-2\lambda _{i}-\left( \lambda
_{i-1}-\lambda _{i}\right) $ and $-2\lambda _{i-1}=-2\lambda _{i}-2\left(
\lambda _{i-1}-\lambda _{i}\right) $ belong to $\Pi _{\Theta }\left(
-2\lambda _{i}\right) $. Proceeding successively it follows that $-\lambda
_{k}-\lambda _{k+1},-2\lambda _{k}\in \Pi _{\Theta }\left( -2\lambda
_{i}\right) $ if $k=j\left( i\right) ,\ldots ,i-1$. The roots $-\lambda
_{k}-\lambda _{r}$, $j\left( i\right) \leq k<r-1\leq i-1$, also belong to $%
\Pi _{\Theta }\left( -2\lambda _{i}\right) $, since $-\lambda _{k}-\lambda
_{r}=\left( -\lambda _{k}-\lambda _{k+1}\right) +\left( \lambda
_{k+1}-\lambda _{r}\right) $ and $\lambda _{k+1}-\lambda _{r}\in \langle
\Theta \rangle $. Hence $V_{\Theta }\left( -2\lambda _{i}\right) $ contains
the subspace $\sum_{k,r=j\left( i\right) }^{i}\mathfrak{g}_{-\lambda
_{k}-\lambda _{r}}$. This subspace is $\mathfrak{z}_{\Theta }$-invariant
because a root in a connected component of $\Theta $ different from $%
\{\alpha _{j\left( i\right) },\ldots ,\alpha _{i-1}\}$ is orthogonal to the
roots $-\lambda _{k}-\lambda _{r}$, $k,r=j\left( i\right) ,\ldots ,i$.

The last statement is a consequence of the expression for $V_{\Theta }\left(
-2\lambda _{i}\right) $ in (2).
\end{profe}

\vspace{12pt}%

\noindent%
\textbf{Remark:} The first case of the above lemma is included in the second
case by taking $j\left( i\right) =i$.

To look at the representation of $K_{\Theta }$ on the subspace $V\left(
-2\lambda _{i}\right) $ of (\ref{forcelelong}) we make use of the following
geometric meaning: \ Let $\mathfrak{sp}\left( j\left( i\right) ,i\right) $
be the subalgebra generated by the root spaces $\mathfrak{g}_{\pm \left(
\lambda _{k}\pm \lambda _{r}\right) }$, $k,r=j\left( i\right) ,\ldots ,i$.
Its elements are symplectic matrices 
\begin{equation*}
\left( 
\begin{array}{cc}
A & -B \\ 
B & A%
\end{array}%
\right)
\end{equation*}%
with $A$, $B$ and $C$ having non zero entries only at the positions $%
k,r=j\left( i\right) ,\ldots ,i$, which shows that it is isomorphic to the
Lie algebra of symplectic matrices in the subspace spanned by $\{e_{j\left(
i\right) },\ldots ,e_{i},f_{j\left( i\right) },\ldots ,f_{i}\}$. Let $%
\mathrm{Sp}\left( j\left( i\right) ,i\right) =\langle \exp \mathfrak{sp}%
\left( j\left( i\right) ,i\right) \rangle $ be the connected subgroup with
Lie algebra $\mathfrak{sp}\left( j\left( i\right) ,i\right) $ and put $%
\mathrm{U}\left( j\left( i\right) ,i\right) =\mathrm{Sp}\left( j\left(
i\right) ,i\right) \cap K$ for its maximal compact subgroup, which is
isomorphic to the unitarian group $\mathrm{U}\left( i-j\left( i\right)
+1\right) $.

The inclusion $\{\alpha _{j\left( i\right) },\ldots ,\alpha _{i-1}\}\subset
\Theta $ shows that the root spaces $\mathfrak{g}_{\lambda _{k}-\lambda
_{r}} $, $k,r=j\left( i\right) ,\ldots ,i$, are contained in the isotropy
subalgebra at the origin $b_{\Theta }\in \mathbb{F}_{\Theta }$. From this it
is \ easily seen that the orbit $\mathrm{Sp}\left( j\left( i\right)
,i\right) \cdot b_{\Theta }=\mathrm{U}\left( j\left( i\right) ,i\right)
\cdot b_{\Theta }$ is a flag manifold of $\mathrm{Sp}\left( j\left( i\right)
,i\right) $ and identifies to the coset $\mathrm{U}\left( j\left( i\right)
,i\right) /\mathrm{SO}\left( j\left( i\right) ,i\right) $ where $\mathrm{SO}%
\left( j\left( i\right) ,i\right) $ is the subgroup isomorphic to $\mathrm{SO%
}\left( i-j\left( i\right) +1\right) $, whose Lie algebra is contained in $%
\sum_{k,r=j\left( i\right) }^{i}\mathfrak{g}_{\lambda _{k}-\lambda _{r}}$.
Further $V_{\Theta }\left( -2\lambda _{i}\right) =\sum_{k,r=j\left( i\right)
}^{i}\mathfrak{g}_{-\lambda _{k}-\lambda _{r}}$ is the tangent space at the
origin $b_{\Theta }\in \mathbb{F}_{\Theta }$ of the orbit $\mathrm{Sp}\left(
j\left( i\right) ,i\right) \cdot b_{\Theta }=\mathrm{U}\left( j\left(
i\right) ,i\right) \cdot b_{\Theta }$.

Now we can get the $K_{\Theta }$-irreducible components of $V_{\Theta
}\left( -2\lambda _{i}\right) $.

\begin{lema}
\label{lemceledecomp}The $\mathfrak{z}_{\Theta }$-irreducible subspace $%
V_{\Theta }\left( -2\lambda _{i}\right) =\sum_{k,r=j\left( i\right) }^{i}%
\mathfrak{g}_{-\lambda _{k}-\lambda _{r}}$ has two $K_{\Theta }$-irreducible
components if $j\left( i\right) <i$. They are given as follows:

\begin{enumerate}
\item The one dimensional subspace $V_{\Theta }\left( -2\lambda _{i}\right)
_{\mathrm{cent}}\subset \sum_{k=j\left( i\right) }^{i}\mathfrak{g}%
_{-2\lambda _{k}}$ spanned by the matrix 
\begin{equation*}
\left( 
\begin{array}{cc}
0 & 0 \\ 
I_{\left[ j\left( i\right) ,i\right] } & 0%
\end{array}%
\right) \in \mathfrak{sp}\left( j\left( i\right) ,i\right) \subset \mathfrak{%
sp}\left( l,\mathbb{R}\right)
\end{equation*}%
where $I_{\left( j\left( i\right) ,i\right) }$ is the diagonal matrix with $%
1 $ in the positions $j\left( i\right) ,\ldots ,i$ and $0$ otherwise.

\item The subspace $V_{\Theta }\left( -2\lambda _{i}\right) _{\mathrm{su}%
\left( j\left( i\right) ,i\right) }$ given by the matrices 
\begin{equation*}
\left( 
\begin{array}{cc}
A & 0 \\ 
B & -A^{T}%
\end{array}%
\right) \in \mathfrak{sp}\left( j\left( i\right) ,i\right)
\end{equation*}%
with $A$ lower triangular and $\mathrm{tr}B=0$.
\end{enumerate}
\end{lema}

\begin{profe}
The compact group $\mathrm{U}\left( j\left( i\right) ,i\right) $ being
isomorphic to $\mathrm{U}\left( i-j\left( i\right) +1\right) $ is the
product of its center $Z_{\left( j\left( i\right) ,i\right) }$ by $\mathrm{SU%
}\left( j\left( i\right) ,i\right) $. The Lie algebra of the center is given
by matrices 
\begin{equation*}
\left( 
\begin{array}{cc}
0 & -B \\ 
B & 0%
\end{array}%
\right)
\end{equation*}%
with $B\in \mathbb{R}\cdot I_{\left( \left( j\left( i\right) ,i\right)
\right) }$ (corresponding to the scalar matrices in $\mathfrak{u}\left(
i-j\left( i\right) +1\right) $). The Lie algebra $\mathfrak{su}\left( \left(
j\left( i\right) ,i\right) \right) $ of $\mathrm{SU}\left( j\left( i\right)
,i\right) $ is given by matrices 
\begin{equation*}
\left( 
\begin{array}{cc}
A & -B \\ 
B & A%
\end{array}%
\right) \in \mathfrak{sp}\left( j\left( i\right) ,i\right)
\end{equation*}%
with $A$ skew symmetric and $B$ symmetric with $\mathrm{tr}B=0$. It follows
that the tangent spaces to the orbits $Z_{\left( j\left( i\right) ,i\right)
}\cdot b_{\Theta }$ and $\mathrm{SU}\left( j\left( i\right) ,i\right) $ are $%
V_{\Theta }\left( -2\lambda _{i}\right) _{\mathrm{cent}}$ and $V_{\Theta
}\left( -2\lambda _{i}\right) _{\mathrm{su}\left( j\left( i\right) ,i\right)
}$, respectively.

Hence, by Lemma \ref{lemorbitnormal} these subspaces are invariant by the
isotropy representation of $\mathrm{SO}\left( j\left( i\right) ,i\right) =%
\mathrm{U}\left( j\left( i\right) ,i\right) \cap K_{\Theta }$. They are $%
K_{\Theta }$-invariant as well because the connected components of $\Theta $
besides $\{\alpha _{j\left( i\right) },\ldots ,\alpha _{i-1}\}$ are
orthogonal to $\Pi _{\Theta }\left( -2\lambda _{i}\right) $. Hence the
simple components of $K_{\Theta }$ different from $\mathrm{SO}\left( j\left(
i\right) ,i\right) $ act trivially on $V_{\Theta }\left( -2\lambda
_{i}\right) $.

Finally, both subspaces $V_{\Theta }\left( -2\lambda _{i}\right) _{\mathrm{%
cent}}$ and $V_{\Theta }\left( -2\lambda _{i}\right) _{\mathrm{su}\left(
j\left( i\right) ,i\right) }$ are irreducible. This is obvious for $%
V_{\Theta }\left( -2\lambda _{i}\right) _{\mathrm{cent}}$ which is one
dimensional. On the other hand the representation of $\mathrm{SO}\left(
j\left( i\right) ,i\right) $ on $V_{\Theta }\left( -2\lambda _{i}\right) $
is equivalent to the isotropy representation of the symmetric space $\mathrm{%
U}\left( j\left( i\right) ,i\right) /\mathrm{SO}\left( j\left( i\right)
,i\right) $, which is known to be irreducible.
\end{profe}

The $\mathfrak{z}_{\Theta }$-irreducible components described in Lemma \ref%
{lemcelelong} contain all the root spaces of the long roots not in $\langle
\Theta \rangle $. They include also the short roots $-\lambda _{i}-\lambda
_{j}$ such that $\lambda _{i}-\lambda _{j}\in \langle \Theta \rangle $. The
other $\mathfrak{z}_{\Theta }$-components are given as follows.

\begin{lema}
\label{lemceleshort}Suppose the root $\lambda _{i}-\lambda _{j}$, $i<j$,
does not belong to $\langle \Theta \rangle $.

\begin{enumerate}
\item If $2\lambda _{l}\notin \Theta $ then the set $\Pi _{\Theta }\left(
-\lambda _{i}+\lambda _{j}\right) $ corresponding to the $\mathfrak{z}%
_{\Theta }$-irreducible component $V_{\Theta }\left( -\lambda _{i}+\lambda
_{j}\right) $ contains only short roots of the type $-\lambda _{r}+\lambda
_{s}$. Similarly $\Pi _{\Theta }\left( -\lambda _{i}-\lambda _{j}\right) $
contains only roots of the type $-\lambda _{r}-\lambda _{s}$.

\item In case $2\lambda _{l}\in \Theta $ let $i_{0}$ be such that $%
C_{l-i_{0}+1}=\{\alpha _{i_{0}+1},\ldots ,\alpha _{l}=2\lambda _{l}\}$ is
the connected component of $\Theta $ containing $2\lambda _{l}$. Then $%
V_{\Theta }\left( -\lambda _{i}+\lambda _{j}\right) =V_{\Theta \setminus
\{2\lambda _{l}\}}\left( -\lambda _{i}+\lambda _{j}\right) $ and $V_{\Theta
}\left( -\lambda _{i}-\lambda _{j}\right) =V_{\Theta \setminus \{2\lambda
_{l}\}}\left( -\lambda _{i}-\lambda _{j}\right) $ if $j\leq i_{0}$.

\item On the other hand if $j\geq i_{0}+1$ then 
\begin{eqnarray}
V_{\Theta }\left( -\lambda _{i}+\lambda _{j}\right) &=&V_{\Theta }\left(
-\lambda _{i}-\lambda _{j}\right)  \label{forceledirectsum} \\
&=&V_{\Theta \setminus C_{l-i_{0}+1}}\left( -\lambda _{i}+\lambda
_{j}\right) \oplus V_{\Theta \setminus C_{l-i_{0}+1}}\left( -\lambda
_{i}-\lambda _{j}\right) .  \notag
\end{eqnarray}
\end{enumerate}

Moreover these $\mathfrak{z}_{\Theta }$-irreducible components are $%
K_{\Theta }$-irreducible.
\end{lema}

\begin{profe}
The first statement is proved as Lemma \ref{lembelecompseml} for $B_{l}$.
The proof of Lemma \ref{lembelecompcoml} also works for the components in
(2) with $j\leq i_{0}$. The direct sum in (\ref{forceledirectsum}) \ comes
from the pair of roots $-\lambda _{i}+\lambda _{j}$, $\left( -\lambda
_{i}+\lambda _{j}\right) -2\lambda _{j}$ and the fact that $2\lambda _{j}\in
\langle \Theta \rangle ^{-}$ if $j\geq i_{0}+1$.

The $K_{\Theta }$-irreducibility of the subspaces in (1) is a consequence of
Lemma \ref{lemirreducible}. In fact no two roots in $\Pi _{\Theta }\left(
-\lambda _{i}+\lambda _{j}\right) $ or in $\Pi _{\Theta }\left( -\lambda
_{i}-\lambda _{j}\right) $ are $M$-equivalent and by Lemma \ref%
{lemtransimplelong} 2b, $\mathcal{W}_{\Theta }$ acts transitively on these
sets of short roots.

The same argument hold for the subspaces $V_{\Theta }\left( -\lambda
_{i}+\lambda _{j}\right) =V_{\Theta \setminus \{2\lambda _{l}\}}\left(
-\lambda _{i}+\lambda _{j}\right) $ and $V_{\Theta }\left( -\lambda
_{i}-\lambda _{j}\right) =V_{\Theta \setminus \{2\lambda _{l}\}}\left(
-\lambda _{i}-\lambda _{j}\right) $ when $j\leq i_{0}$, which are indeed $%
K_{\Theta \setminus \{2\lambda _{l}\}}$-irreducible.

Finally, in (2) if $j\geq i_{0}+1$ then $2\lambda _{j}\in \langle \Theta
\rangle $. From the equality $-\lambda _{i}-\lambda _{j}=-\lambda
_{i}+\lambda _{j}-2\lambda _{j}$ we see that $\mathfrak{g}_{-\lambda
_{i}+\lambda _{l}}\oplus \mathfrak{g}_{-\lambda _{i}-\lambda _{l}}$ is an
irreducible subspace for the three dimensional subalgebra $\mathfrak{g}%
\left( 2\lambda _{j}\right) $, isomorphic to $\mathfrak{sl}\left( 2,\mathbb{R%
}\right) $, generated by $\mathfrak{g}_{\pm 2\lambda _{j}}$. In this two
dimensional subspace the compact part $\mathfrak{k}\left( 2\lambda
_{j}\right) $ of $\mathfrak{g}\left( 2\lambda _{j}\right) $ is also
irreducible. Now let $\{0\}\neq U\subset V_{\Theta }\left( -\lambda
_{i}+\lambda _{j}\right) $ be a $K_{\Theta }$-invariant subspace. Since the
roots $-\lambda _{i}\pm \lambda _{j}\in \Pi _{\Theta }\left( -\lambda
_{i}+\lambda _{j}\right) =\Pi _{\Theta }\left( -\lambda _{i}-\lambda
_{j}\right) $ are $M$-equivalent, it follows by Lemma \ref{leminterinvMclass}
that $U\cap \left( \mathfrak{g}_{-\lambda _{i}+\lambda _{l}}\oplus \mathfrak{%
g}_{-\lambda _{i}-\lambda _{l}}\right) \neq \{0\}$. This subspace is $%
\mathfrak{k}\left( 2\lambda _{j}\right) $-invariant, hence $\mathfrak{g}%
_{-\lambda _{i}+\lambda _{l}}\oplus \mathfrak{g}_{-\lambda _{i}-\lambda
_{l}}\subset U$. Hence irrducibility of $V_{\Theta }\left( -\lambda
_{i}+\lambda _{j}\right) $ is a consequence of Lemma \ref{lemirreducible}
combined with the fact that $\mathcal{W}_{\Theta \setminus \{2\lambda
_{l}\}}\subset \mathcal{W}_{\Theta }$ acts transitively on the sets $\Pi
_{\Theta \setminus \{2\lambda _{l}\}}\left( -\lambda _{i}+\lambda
_{j}\right) $ and $\Pi _{\Theta \setminus \{2\lambda _{l}\}}\left( -\lambda
_{i}-\lambda _{j}\right) $.
\end{profe}

With this lemma we finish the description of the irreducible $K_{\Theta }$%
-components. Among them the only $K_{\Theta }$-equivalents are the following:

\begin{enumerate}
\item The one dimensional subspaces $V_{\Theta }\left( -2\lambda _{i}\right)
_{\mathrm{cent}}$ of Lemma \ref{lemceledecomp} (1). The representation of $%
K_{\Theta }$ on each one of them is trivial.

\item $V_{\Theta }\left( -\lambda _{i}+\lambda _{j}\right) \approx V_{\Theta
}\left( -\lambda _{i}-\lambda _{j}\right) $ when $2\lambda _{l}\notin \Theta 
$ as in Lemma \ref{lemceleshort} (1). This equivalence follows by
Proposition \ref{propequivalent}, since there is a bijection between $\Pi
_{\Theta }\left( -\lambda _{i}+\lambda _{j}\right) $ and $\Pi _{\Theta
}\left( -\lambda _{i}-\lambda _{j}\right) $ mapping a root $-\lambda
_{r}+\lambda _{s}\in \Pi _{\Theta }\left( -\lambda _{i}+\lambda _{j}\right) $
to the $M$-equivalent root $-\lambda _{r}-\lambda _{s}\in \Pi _{\Theta
}\left( -\lambda _{i}-\lambda _{j}\right) $.

\item The subspaces $V_{\Theta }\left( -\lambda _{i}+\lambda _{j}\right)
=V_{\Theta \setminus \{2\lambda _{l}\}}\left( -\lambda _{i}+\lambda
_{j}\right) $ and $V_{\Theta }\left( -\lambda _{i}-\lambda _{j}\right)
=V_{\Theta \setminus \{2\lambda _{l}\}}\left( -\lambda _{i}-\lambda
_{j}\right) $ with $j\leq i_{0}$ as in Lemma \ref{lemceleshort} (2).
\end{enumerate}

Any other pair of subspaces are not $K_{\Theta }$-equivalent because the
lack of $M$-equivalence in the corresponding sets of roots (cf. Lemma \ref%
{lemnotequivalent}).

Summarizing we get the $K_{\Theta }$-invariant subspaces for the flags of $%
C_{l}$, $l>4$.

\begin{teorema}
Let $\mathbb{F}_{\Theta }$ be a flag manifold of $C_{l}=\mathfrak{sp}\left(
l,\mathbb{R}\right) $, $l\geq 5$. Then the $K_{\Theta }$-invariant
irreducible subspaces of $\mathfrak{n}_{\Theta }^{-}$ are the following:

\begin{enumerate}
\item Continuous families:

\begin{enumerate}
\item One dimensional subspaces spanned by matrices 
\begin{equation*}
\left( 
\begin{array}{cc}
0 & 0 \\ 
\Lambda & 0%
\end{array}%
\right)
\end{equation*}%
where $\Lambda $ is a diagonal matrix $a_{1}I_{\left[ j\left( i_{1}\right)
,i_{1}\right] }+\cdots +a_{k}I_{\left[ j\left( i_{k}\right) ,i_{k}\right] }$
where $\left[ j\left( i_{1}\right) ,i_{1}\right] $, \ldots ,$\left[ j\left(
i_{1}\right) ,i_{1}\right] $ are the connected components of $\Theta $ not
containing $2\lambda _{l}$, and $I_{\left[ j\left( i_{s}\right) ,i_{s}\right]
}$ is the identity matrix corresponing to these indexes.

\item The subspaces parametrized by $[(x,y)]\in \mathbb{R}P^{2}$ given by 
\begin{equation*}
V_{[(x,y)]}^{ij}=\left\{ xX+yTX:X\in V_{\Theta }\left( -\lambda _{i}+\lambda
_{j}\right) \right\} \qquad i<j
\end{equation*}%
where $T:V_{\Theta }\left( -\lambda _{i}+\lambda _{j}\right) \rightarrow
V_{\Theta }\left( -\lambda _{i}-\lambda _{j}\right) $ is an itertiwining
operator for the $K_{\Theta }$-representations. Here the indexes $ij$ are
arbitrary if $2\lambda _{l}\notin \Theta $. Otherwise $j\leq i_{0}$, where $%
\{\alpha _{i_{0}+1},\ldots ,\alpha _{l}=2\lambda _{l}\}$ is the component of 
$\Theta $ containing $2\lambda _{l}$.
\end{enumerate}

\item Isolated subspaces:

\begin{enumerate}
\item The subspaces $V_{\Theta }\left( -2\lambda _{i}\right) _{\mathrm{su}%
\left( j\left( i\right) ,i\right) }$ of codimension $1$ contained in $%
V_{\Theta }\left( -2\lambda _{i}\right) $ as defined in Lemma \ref%
{lemceledecomp}.

\item The subspaces%
\begin{eqnarray*}
V_{\Theta }\left( -\lambda _{i}+\lambda _{j}\right) &=&V_{\Theta }\left(
-\lambda _{i}-\lambda _{j}\right) \\
&=&V_{\Theta \setminus C_{l-i_{0}+1}}\left( -\lambda _{i}+\lambda
_{j}\right) \oplus V_{\Theta \setminus C_{l-i_{0}+1}}\left( -\lambda
_{i}-\lambda _{j}\right) ,
\end{eqnarray*}%
when $2\lambda _{l}\in \Theta $ and $i<i_{0}+1\leq j$ where $\{\alpha
_{i_{0}+1},\ldots ,\alpha _{l}=2\lambda _{l}\}$ is the component of $\Theta $
containing $2\lambda _{l}$.
\end{enumerate}
\end{enumerate}
\end{teorema}

When $l=4$ the $M$-equivalence classes of the short roots increase to $%
\{\lambda _{1}-\lambda _{2},\lambda _{1}+\lambda _{2},\lambda _{3}-\lambda
_{4},\lambda _{3}+\lambda _{4}\}$, $\{\lambda _{1}-\lambda _{3},\lambda
_{1}+\lambda _{3},\lambda _{2}-\lambda _{4},\lambda _{2}+\lambda _{4}\}$, $%
\{\lambda _{1}-\lambda _{4},\lambda _{1}+\lambda _{4},\lambda _{2}-\lambda
_{3},\lambda _{2}+\lambda _{3}\}$ while the long roots are kept the same $%
\{2\lambda _{1},2\lambda _{2},2\lambda _{3},2\lambda _{4}\}$. Since there
are more $M$-equivalent pair of roots we can have more $K_{\Theta }$%
-equivalent subspaces than in the general case.

For example consider flag $\mathbb{F}_{\{\lambda _{2}-\lambda _{3}\}}$. By
the general result the subspaces $V_{\{\lambda _{2}-\lambda _{3}\}}\left(
-\lambda _{1}+\lambda _{2}\right) $ and $V_{\{\lambda _{2}-\lambda
_{3}\}}\left( -\lambda _{1}-\lambda _{2}\right) $ are $K_{\Theta }$%
-equivalent. Their corresponding roots are $\Pi _{\{\lambda _{2}-\lambda
_{3}\}}\left( -\lambda _{1}+\lambda _{2}\right) =\{-\lambda _{1}+\lambda
_{2},-\lambda _{1}+\lambda _{3}\}$ and $\Pi _{\{\lambda _{2}-\lambda
_{3}\}}\left( -\lambda _{1}-\lambda _{2}\right) =\{-\lambda _{1}-\lambda
_{2},-\lambda _{1}-\lambda _{3}\}$. For $l=4$ we have $\left( -\lambda
_{1}+\lambda _{2}\right) \sim _{M}\left( -\lambda _{3}+\lambda _{4}\right) $
and $\left( -\lambda _{1}+\lambda _{3}\right) \sim _{M}\left( -\lambda
_{2}+\lambda _{4}\right) $. Since the set of root for $V_{\{\lambda
_{2}-\lambda _{3}\}}\left( -\lambda _{3}+\lambda _{4}\right) $ is $\Pi
_{\{\lambda _{2}-\lambda _{3}\}}\left( -\lambda _{3}+\lambda _{4}\right)
=\{-\lambda _{3}+\lambda _{4},-\lambda _{2}+\lambda _{4}\}$ we conclude that 
$V_{\{\lambda _{2}-\lambda _{3}\}}\left( -\lambda _{1}+\lambda _{2}\right) $
is also $K_{\Theta }$-equivalent to $V_{\{\lambda _{2}-\lambda _{3}\}}\left(
-\lambda _{3}+\lambda _{4}\right) $. The same way $V_{\{\lambda _{2}-\lambda
_{3}\}}\left( -\lambda _{1}-\lambda _{2}\right) $ and $V_{\{\lambda
_{2}-\lambda _{3}\}}\left( -\lambda _{3}-\lambda _{4}\right) $ are $%
K_{\Theta }$-equivalent. Thus we get new continuous families of invariant
subspaces that are not present in the general case.

\subsection{Flags of $D_{l}=\mathfrak{so}\left( l,l\right) $}

The Dynkin diagram $D_{l}$ has no multiple edges. Hence on any flag manifold 
$\mathbb{F}_{\Theta }$ we have, by Lemma \ref{lemtransimplelong}, that $%
\mathcal{W}_{\Theta }$ acst transitively on each set of roots $\Pi _{\Theta
}^{\sigma }$ corresponding to an irreducible representation of $\mathfrak{z}%
_{\Theta }$ on $V_{\Theta }^{\sigma }\subset \mathfrak{n}_{\Theta }^{-}$.
This transitivity is one of the conditions of Lemma \ref{lemirreducible},
ensuring that the subspaces $V_{\Theta }^{\sigma }$ are $K_{\Theta }$%
-irreducible. To look at the condition involving the $M$-equivalence classes
we work, as in Section \ref{secmequiv}, with the standard realization of $%
D_{l}=\mathfrak{so}\left( l,l\right) $.

In this realization $\mathfrak{a}$ is the subalgebra of matrices 
\begin{equation*}
\left( 
\begin{array}{cc}
\Lambda & 0 \\ 
0 & -\Lambda%
\end{array}%
\right)
\end{equation*}%
with $\Lambda =\mathrm{diag}\{a_{1},\ldots ,a_{l}\}$ and the set of simple
roots is $\Sigma =\{\lambda _{1}-\lambda _{2},\ldots ,\lambda _{l-1}-\lambda
_{l},\lambda _{l-1}+\lambda _{l}\}$.

The Weyl chamber $\mathfrak{a}^{+}\subset \mathfrak{a}$ is defined by the
inequalities 
\begin{equation}
a_{1}>a_{2}>\cdots >a_{l-1}>a_{l}>-a_{l-1}.  \label{forinechamber}
\end{equation}%
A partial chamber $\mathfrak{a}_{\Theta }\cap \mathrm{cl}\mathfrak{a}^{+}$
is defined by a similar relation where some of the inequalities are changed
by equalities. In particular a characteristic element $H_{\Theta }$ for the
subset $\Theta =\{\alpha \in \Sigma :\alpha \left( H_{\Theta }\right)
=0\}\subset \Sigma $ is defined by one of these relations.

The following statement is specific for $D_{l}$ and will be used soon to
check that $M$-equivalent root spaces are not contained in an irreducible
component.

\begin{lema}
Given a subset $\Theta \subset \Sigma $ there exists characteristic element 
\begin{equation*}
H_{\Theta }=\left( 
\begin{array}{cc}
\Lambda _{\Theta } & 0 \\ 
0 & -\Lambda _{\Theta }%
\end{array}%
\right) \in \mathfrak{a}_{\Theta }\cap \mathrm{cl}\mathfrak{a}^{+}
\end{equation*}%
with $\Lambda _{\Theta }=\mathrm{diag}\{a_{1},\ldots ,a_{l}\}$ such that $%
a_{i}\neq 0$, $i=1,\ldots ,l$.
\end{lema}

\begin{profe}
By the last two inequalities in (\ref{forinechamber}) we have $a_{l-1}\geq
-a_{l-1}$, that is, $a_{l-1}\geq 0$. Also, $a_{l-1}=0$ if and only if $%
a_{l-1}=a_{l}=-a_{l-1}$, that is, $a_{l-1}-a_{l}=a_{l-1}+a_{l}=0$, which
means that both roots $\lambda _{l-1}-\lambda _{l}$ and $\lambda
_{l-1}+\lambda _{l}$ belong to $\Theta $. This being so we consider the
possibilities:

\begin{enumerate}
\item $\{\lambda _{l-1}-\lambda _{l},\lambda _{l-1}+\lambda _{l}\}\subset
\Theta $. Let $i<l-1$ be the maximum such that $\lambda _{i}-\lambda
_{i+1}\notin \Theta $ (it is tacitly assumed that $\Theta \neq \Sigma $).
Then the conditions to define a characteristic element for $\Theta $ have
the form%
\begin{equation*}
a_{1}\geq \cdots \geq a_{i}>a_{i+1}=a_{i+2}=\cdots =a_{l}.
\end{equation*}%
Thus we can choose a characteristic element having $a_{i+1}=a_{i+2}=\cdots
=a_{l}>0$, so that all the entries of $\Lambda _{\Theta }$ will be $>0$.

\item One of the roots $\lambda _{l-1}-\lambda _{l}$ or $\lambda
_{l-1}+\lambda _{l}$ does not belong to $\Theta $. In this case $%
a_{l-1}>0>-a_{l-1}$ for any $H_{\Theta }$ so that $a_{i}>0$ for any $i\leq
l-1$. Also the relations defining $\mathfrak{a}_{\Theta }\cap \mathrm{cl}%
\mathfrak{a}^{+}$ end with 
\begin{equation*}
a_{l-1}>a_{l}>-a_{l-1}\quad \mathrm{or}\quad a_{l-1}\geq a_{l}>-a_{l-1}\quad 
\mathrm{or}\quad a_{l-1}>a_{l}\geq -a_{l-1}.
\end{equation*}%
In each case we can choose $a_{l}\neq 0$ without violating the conditions.
\end{enumerate}
\end{profe}

From now on we distinguish the cases where $l>4$ and $l=4$.

If $l>4$ then $M$-equivalence classes in the positive roots are $\{\lambda
_{i}-\lambda _{j},\lambda _{i}+\lambda _{j}\}$, $i<j$, and $\{\lambda
_{i}-\lambda _{j},-\lambda _{i}-\lambda _{j}\}$, $i>j$, in the negative
roots. By the previous lemma we get easily that the corresponding $M$%
-equivalent root spaces are not contained in a single irreducible component.

\begin{lema}
\label{lemdelecompequiv}Let $V_{\Theta }^{\sigma }$ be an irreducible
component containing the root spaces $\mathfrak{g}_{\alpha }$ and $\mathfrak{%
g}_{\beta }$, $\alpha \neq \beta $. If $l>4$ then $\alpha $ and $\beta $ are
not $M$-equivalent.
\end{lema}

\begin{profe}
Take a characteristic element $H_{\Theta }$ with $a_{i}\neq 0$ as in the
previous lemma. Then $\left( \lambda _{i}-\lambda _{j}\right) \left(
H_{\Theta }\right) \neq -\left( \lambda _{i}+\lambda _{j}\right) \left(
H_{\Theta }\right) $ for otherwise $a_{i}-a_{j}=-a_{i}-a_{j}$, that is, $%
a_{i}=0$. The result follows, since $V_{\Theta }^{\sigma }$ is contained in
an eigenspace of $\mathrm{ad}\left( H_{\Theta }\right) $.
\end{profe}

Combining this lemma with Lemma \ref{lemtransimplelong} (about the
transitivity of $\mathcal{W}_{\Theta }$) we get at once $K_{\Theta }$%
-irreducibility of $V_{\Theta }^{\sigma }$, by Lemma \ref{lemirreducible}.

\begin{proposicao}
In any flag manifold $\mathbb{F}_{\Theta }$ of $D_{l}$, $l>4$, the $%
\mathfrak{z}_{\Theta }$-irreducible components $V_{\Theta }^{\sigma }$ are
also $K_{\Theta }$-irreducible.
\end{proposicao}

To get the full picture of the invariant subspaces we must find the pairs of 
$\mathfrak{z}_{\Theta }$-irreducible components that are mutually $K_{\Theta
}$-equivalent. Our method to check $K_{\Theta }$-equivalence is via the
orbits on $\mathbb{F}_{\Theta }$ of the simple components of the maximal
compact subgroup $K$ of $G$.

To this purpose we need some further notation concerning the standard
realization of $D_{l}$. The Lie algebra $\mathfrak{so}\left( l,l\right) $ is
the algebra of $2l\times 2l$ matrices of the form 
\begin{equation*}
\left( 
\begin{array}{cc}
A & B \\ 
C & -A^{T}%
\end{array}%
\right) \quad B+B^{T}=C+C^{T}=0.
\end{equation*}%
We have that $\mathfrak{k}$ is the subalgebra of skew-symmetric matrices in $%
\mathfrak{so}\left( l,l\right) $, that is,%
\begin{equation*}
\left( 
\begin{array}{cc}
A & B \\ 
B & A%
\end{array}%
\right) \quad A+A^{T}=B+B^{T}=0.
\end{equation*}%
$\mathfrak{k}$ is the direct sum of two copies of $\mathfrak{so}\left(
l\right) $. In fact, via the decomposition 
\begin{equation*}
\left( 
\begin{array}{cc}
A & B \\ 
B & A%
\end{array}%
\right) =\left( 
\begin{array}{cc}
\left( A+B\right) /2 & \left( A+B\right) /2 \\ 
\left( A+B\right) /2 & \left( A+B\right) /2%
\end{array}%
\right) +\left( 
\begin{array}{cc}
\left( A-B\right) /2 & -\left( A-B\right) /2 \\ 
-\left( A-B\right) /2 & \left( A-B\right) /2%
\end{array}%
\right)
\end{equation*}%
we get $\mathfrak{k}=\mathfrak{so}\left( l\right) _{1}\oplus \mathfrak{so}%
\left( l\right) _{2}$ with%
\begin{equation*}
\mathfrak{so}\left( l\right) _{1}:\left( 
\begin{array}{ll}
A & A \\ 
A & A%
\end{array}%
\right) \qquad \mathfrak{so}\left( l\right) _{2}:\left( 
\begin{array}{cc}
A & -A \\ 
-A & A%
\end{array}%
\right)
\end{equation*}%
where in both cases $A$ is skew-symmetric. We write $\mathrm{SO}\left(
l\right) _{i}=\langle \exp \mathfrak{so}\left( l\right) _{i}\rangle $, $%
i=1,2 $. \ 

As to the root spaces we write 
\begin{equation}
E_{ij}^{-}=\left( 
\begin{array}{cc}
E_{ij} & 0 \\ 
0 & -E_{ij}^{t}%
\end{array}%
\right) \qquad \mathrm{and}\qquad E_{ij}^{+}=\left( 
\begin{array}{cc}
0 & 0 \\ 
E_{ij}-E_{ij}^{t} & 0%
\end{array}%
\right)  \label{forrootdele}
\end{equation}%
where $E_{ij}$ is a basic $l\times l$ matrix. Then $E_{ij}^{-}$ spans the
root space $\mathfrak{g}_{\lambda _{i}-\lambda _{j}}$ and $E_{ij}^{+}$ spans 
$\mathfrak{g}_{-\lambda _{i}-\lambda _{j}}$.

We can return now to the question of $K_{\Theta }$-equivalence of the $%
\mathfrak{z}_{\Theta }$-components $V_{\Theta }^{\sigma }$.

For a root $\alpha \in \Pi ^{-}\setminus \langle \Theta \rangle ^{-}$ write $%
V_{\Theta }\left( \alpha \right) $ for the irreducible component containing $%
\mathfrak{g}_{\alpha }$ (cf. Proposition \ref{propcomproot}). By Lemma \ref%
{lemdelecompequiv} we have $V_{\Theta }\left( \alpha \right) \neq V_{\Theta
}\left( \beta \right) $ if $\alpha \sim _{M}\beta $ and $\alpha \neq \beta $%
. Moreover, by Lemma \ref{lemnotequivalent} a component $V_{\Theta }^{\sigma
}$ is not $K_{\Theta }$-equivalent to $V_{\Theta }\left( \alpha \right) $
unless there exists $\beta \sim _{M}\alpha $ such that $V_{\Theta }^{\sigma
}=V_{\Theta }\left( \beta \right) $.

Now by Section \ref{secmequiv} we have that if $l\neq 4$ then the $M$%
-equivalent classes of $D_{l}$ (on the negative roots) has exactly two
elements. If $\{\alpha ,\beta \}$ is a $M$-equivalence class with say $%
\alpha \in \Pi ^{-}\setminus \langle \Theta \rangle ^{-}$ and $\beta \in
\langle \Theta \rangle ^{-}$ then $V_{\Theta }\left( \alpha \right) $ is not 
$K_{\Theta }$-equivalent to any other irreducible component $V_{\Theta
}^{\sigma }$. On the other hand if both $\alpha ,\beta \in \Pi ^{-}\setminus
\langle \Theta \rangle ^{-}$ we have $K_{\Theta }$-equivalence between $%
V_{\Theta }\left( \alpha \right) $ and $V_{\Theta }\left( \beta \right) $.

\begin{lema}
In $D_{l}$, $l>4$, let $\{\alpha ,\beta \}$ be a $M$-equivalence class
contained in $\Pi ^{-}\setminus \langle \Theta \rangle ^{-}$. Then the $%
K_{\Theta }$ representations on $V_{\Theta }\left( \alpha \right) $ and $%
V_{\Theta }\left( \beta \right) $ are equivalent.
\end{lema}

\begin{profe}
To prove equivalence we shall exhibit a $K_{\Theta }$-invariant subspace $%
\{0\}\neq V\subset V_{\Theta }\left( \alpha \right) \oplus V_{\Theta }\left(
\beta \right) $ which is different from the irreducible components $%
V_{\Theta }\left( \alpha \right) $ and $V_{\Theta }\left( \beta \right) $.
This will imply that the components are indeed $K_{\Theta }$-equivalent (see
Proposition \ref{propsumirredcomp}).

The required subspace $V$ will be obtained from the tangent space at the
origin of the orbit of one of the normal subgroups $\mathrm{SO}\left(
l\right) _{j}$, $j=1,2$.

Take for instance $\mathrm{SO}\left( l\right) _{1}$ whose Lie algebra $%
\mathfrak{so}\left( l\right) _{1}$ constitutes of the matrices 
\begin{equation*}
\left( 
\begin{array}{cc}
A & A \\ 
A & A%
\end{array}%
\right) \quad A+A^{T}=0.
\end{equation*}%
Looking at these matrices we see that after identifying $T_{b_{\Theta }}%
\mathbb{F}_{\Theta }$ with $\mathfrak{n}_{\Theta }^{-}$ the tangent space $%
T_{b_{\Theta }}\left( \mathrm{SO}\left( l\right) _{1}\cdot b_{\Theta
}\right) $ to the orbit $\mathrm{SO}\left( l\right) _{1}\cdot b_{\Theta }$
is identified to the subspace $W_{1}\subset \mathfrak{n}_{\Theta }^{-}$
spanned by $\mathrm{pr}\left( E_{rs}^{-}+E_{rs}^{+}\right) $, $r>s$, where $%
E_{rs}^{\pm }$ were defined in (\ref{forrootdele}) and $\mathrm{pr}:%
\mathfrak{n}^{-}\rightarrow \mathfrak{n}_{\Theta }^{-}$ is the projection
w.r.t. the root spaces decomposition.

The tangent space $T_{b_{\Theta }}\left( \mathrm{SO}\left( l\right)
_{1}\cdot b_{\Theta }\right) $ is invariant by the isotropy representation
of $K_{\Theta }$, by Lemma \ref{lemorbitnormal}. Hence $W_{1}$ is invariant
by the adjoint action of $K_{\Theta }$.

Now a $M$-equivalence class is given by $\{\lambda _{i}-\lambda
_{j},-\lambda _{i}-\lambda _{j}\}$, $i>j$, whose root spaces are spanned by $%
E_{ij}^{-}$ and $E_{ij}^{+}$, so that $E_{ij}^{-}\in V_{\Theta }\left(
\lambda _{i}-\lambda _{j}\right) $ and $E_{ij}^{+}\in V_{\Theta }\left(
-\lambda _{i}-\lambda _{j}\right) $. If both roots are in $\Pi ^{-}\setminus
\langle \Theta \rangle ^{-}$ we have 
\begin{equation*}
E_{rs}^{-}+E_{rs}^{+}=\mathrm{pr}\left( E_{rs}^{-}+E_{rs}^{+}\right) \in
W_{1}\cap \left( V_{\Theta }\left( \lambda _{i}-\lambda _{j}\right) \oplus
V_{\Theta }\left( -\lambda _{i}-\lambda _{j}\right) \right) .
\end{equation*}%
Hence $W_{1}\cap \left( V_{\Theta }\left( \lambda _{i}-\lambda _{j}\right)
\oplus V_{\Theta }\left( -\lambda _{i}-\lambda _{j}\right) \right) \neq
\{0\} $. This is a $K_{\Theta }$-invariant subspace different from $%
V_{\Theta }\left( \lambda _{i}-\lambda _{j}\right) $ and $V_{\Theta }\left(
-\lambda _{i}-\lambda _{j}\right) $. It follows that the representation of $%
K_{\Theta }$ on the irreducible subspaces $V_{\Theta }\left( \lambda
_{i}-\lambda _{j}\right) $ and $V_{\Theta }\left( -\lambda _{i}-\lambda
_{j}\right) $ are equivalent by Proposition \ref{propsumirredcomp}.
\end{profe}

Summarizing we get the $K_{\Theta }$-invariant subspaces for the flags of $%
D_{l}$, $l>4$.

\begin{teorema}
In a flag $\mathbb{F}_{\Theta }$ of $D_{l}$, $l>4$, there are the following
two classes of $K_{\Theta }$-invariant subspaces in $\mathfrak{n}_{\Theta
}^{-}$.

\begin{enumerate}
\item The $\mathfrak{z}_{\Theta }$-irreducible component $V_{\Theta }\left(
\alpha \right) $, containing the root space $\mathfrak{g}_{\alpha }$ in case 
$\alpha \in \Pi ^{-}\setminus \langle \Theta \rangle ^{-}$ is not $M$%
-equivalent to $\beta \in \Pi ^{-}\setminus \langle \Theta \rangle ^{-}$.
(These are isolated invariant subspaces.)

\item Let $\{\alpha ,\beta \}$ be a $M$-equivalence class contained in $\Pi
^{-}\setminus \langle \Theta \rangle ^{-}$. Then there is a continuum of
invariant subspaces parametrized by $[(x,y)]\in \mathbb{R}P^{1}$ given by 
\begin{equation*}
V_{[(x,y)]}=\left\{ xX+yTX:X\in V_{\Theta }\left( \alpha \right) \right\}
\end{equation*}%
where $T:V_{\Theta }\left( \alpha \right) \rightarrow V_{\Theta }\left(
\beta \right) $ is an isomorphism intertwining the $K_{\Theta }$%
-representations.
\end{enumerate}
\end{teorema}

The case $l=4$ differs from the general one in two aspects, \ namely each $M$%
-equivalence class has now $4$ elements and the compact subalgebra $%
\mathfrak{k}=\mathfrak{so}\left( 4\right) \oplus \mathfrak{so}\left(
4\right) $ decomposes further into four copies of $\mathfrak{so}\left(
3\right) $. These simple components of $\mathfrak{k}$ yield new invariant
subspaces.

To see what can happen let us consider the example with $\Theta =\{\lambda
_{1}-\lambda _{2},\lambda _{2}-\lambda _{3},\lambda _{3}-\lambda _{4}\}$.
Then $\mathfrak{n}_{\Theta }^{-}$ is formed by matrices 
\begin{equation*}
Y=\left( 
\begin{array}{ll}
0 & 0 \\ 
X & 0%
\end{array}%
\right)
\end{equation*}%
with $X$ a $4\times 4$ skew-symmetric matrix. Here $\mathfrak{k}$ is the Lie
algebra of matrices 
\begin{equation*}
\left( 
\begin{array}{cc}
\alpha & \beta \\ 
\beta & \alpha%
\end{array}%
\right) \qquad \alpha +\alpha ^{T}=\beta +\beta ^{T}=0
\end{equation*}%
which is isomorphic to $\mathfrak{so}\left( 4\right) \oplus \mathfrak{so}%
\left( 4\right) $ by 
\begin{equation*}
\left( 
\begin{array}{cc}
\alpha & \beta \\ 
\beta & \alpha%
\end{array}%
\right) \longmapsto \left( 
\begin{array}{cc}
\alpha +\beta & 0 \\ 
0 & \alpha -\beta%
\end{array}%
\right) .
\end{equation*}%
Now $\mathfrak{so}\left( 4\right) $ is the direct sum of two copies of $%
\mathfrak{so}\left( 3\right) $ as in (\ref{fordecomso4so3so3}) (see the case 
$A_{3}$, above). Thus we can see that if we take $X$ in each one of the sets
of matrices in (\ref{fordecomso4so3so3}) we get subespaces $V_{1}$ and $%
V_{2} $ of $\mathfrak{n}_{\Theta }^{-}$ that are the tangent spaces to the
orbits of the simple components of $K$. Hence $\mathfrak{n}_{\Theta
}^{-}=V_{1}\oplus V_{2}$ is a decomposition into two $3$-dimensional $%
K_{\Theta }$-invariant subspaces. These two representations are equivalent,
since they are just the adjoint representation of $\mathfrak{so}\left(
3\right) $ on each component.

\subsection{Flags of $E_{6}$, $E_{7}$ and $E_{8}$}

For a flag manifold $\mathbb{F}_{\Theta }$ of one of these exceptional Lie
algebras the $K_{\Theta }$-invariant irreducible subspaces finite and
coincide with the invariant irreducible components $V_{\Theta }^{\sigma
}\subset \mathfrak{n}_{\Theta }^{-}$ for the representation of $\mathfrak{z}%
_{\Theta }$.

This is because the Dynkin diagrams are simply laced. Hence, by Lemma \ref%
{lemtransimplelong} it follows that $\mathcal{W}_{\Theta }$ acts
transitively on each set of roots $\Pi _{\Theta }^{\sigma }$ corresponding
to $V_{\Theta }^{\sigma }$. Also, as checked in Section \ref{secmequiv} the
classes of $M$-equivalence for these Lie algebras are singletons. Hence by
Lemma \ref{lemirreducible} we have $K_{\Theta }$-irreducibility of each $%
V_{\Theta }^{\sigma }$. Furthermore the representations of $K_{\Theta }$ on
different subspaces $V_{\Theta }^{\sigma _{1}}$ and $V_{\Theta }^{\sigma
_{2}}$ are not $M$-equivalent, as follows by combining Lemma \ref%
{lemnotequivalent} and the fact that the $M$-equivalence classes are
singletons.

\subsection{Flags of $G_{2}$}

Let $\alpha _{1}$ and $\alpha _{2}$ be the simple roots of $G_{2}$ with $%
\alpha _{1}$ the long one. There are three flag manifolds, $\mathbb{F}%
_{\emptyset }$, $\mathbb{F}_{\{\alpha _{1}\}}$ and $\mathbb{F}_{\{\alpha
_{2}\}}$. The irreducible components on them are easily obtained by direct
inspection of the positive roots. Recall that the $M$-equivalence classes on
the positive roots are $\{\alpha _{1},\alpha _{1}+2\alpha _{2}\}$, $\{\alpha
_{1}+\alpha _{2},\alpha _{1}+3\alpha _{2}\}$ and $\{\alpha _{2},2\alpha
_{1}+3\alpha _{2}\}$. They are listed below:

\begin{enumerate}
\item In $\mathbb{F}=\mathbb{F}_{\emptyset }$ there are three families of $%
\mathfrak{z}_{\Theta }$ and $K_{\Theta }$-irreducible subspaces,
parametrized by $\mathbb{R}P^{1}$, corresponding to the three $M$%
-equivalence classes on the negative roots.

\item For $\mathbb{F}_{\{\alpha _{1}\}}$ there are three $\mathfrak{z}%
_{\Theta }$-irreducible components corresponding to the sets of roots $%
\{\alpha _{2},\alpha _{1}+\alpha _{2}\}$, $\{\alpha _{1}+2\alpha _{2}\}$ and 
$\{\alpha _{1}+3\alpha _{2},2\alpha _{1}+3\alpha _{2}\}$. They are $%
K_{\Theta }$-irreducible because the $2$-dimensional irreducible
representation of $\mathfrak{sl}\left( 2,\mathbb{R}\right) $ is $\mathfrak{so%
}\left( 2\right) $-irreducible. By checking the $M$-equivalence classes we
see that the $2$-dimensional subspaces are equivalent. Hence, we have the
irreducible subspace $\mathfrak{g}_{-\alpha _{1}-2\alpha _{2}}$ and a family
of $2$-dimensional irreducible subspaces parametrized by $\mathbb{R}P^{1}$,
contained in $\mathfrak{g}_{-\alpha _{2}}\oplus \mathfrak{g}_{-\alpha
_{1}-\alpha _{2}}\mathfrak{g}_{-\alpha _{1}-3\alpha _{2}}\oplus \mathfrak{g}%
_{-2\alpha _{1}-3\alpha _{2}}$.

\item For $\mathbb{F}_{\{\alpha _{2}\}}$ the $\mathfrak{z}_{\Theta }$%
-irreducible components correspond to the sets of roots $\{\alpha
_{1},\alpha _{1}+\alpha _{2},\alpha _{1}+2\alpha _{2},\alpha _{1}+3\alpha
_{2}\}$ and $\{2\alpha _{1}+3\alpha _{2}\}$. The $4$-dimensional irreducible
representation of $\mathfrak{z}_{\Theta }\approx \mathfrak{sl}\left( 2,%
\mathbb{R}\right) $ decomposes into two $K_{\Theta }\approx \mathrm{SO}%
\left( 2\right) $-invariant irreducible $2$-dimensional inequivalent
representations. Hence there are three $K_{\Theta }$-invariant irreducible
subspaces. \ 
\end{enumerate}

\subsection{Flags of $F_{4}$}

Recall that the $M$-equivalence classes on the positive roots of $F_{4}$ are
given by

\begin{itemize}
\item $12$ singletons $\{\alpha \}$ with $\alpha $ running through the set
of short roots.

\item $3$ sets of long roots $\{2\alpha _{1}+3\alpha _{2}+4\alpha
_{3}+2\alpha _{4},\alpha _{2},\alpha _{2}+2\alpha _{3},\alpha _{2}+2\alpha
_{3}+2\alpha _{4}\}$, $\{\alpha _{1}+3\alpha _{2}+4\alpha _{3}+2\alpha
_{4},\alpha _{1}+\alpha _{2},\alpha _{1}+\alpha _{2}+2\alpha _{3},\alpha
_{1}+\alpha _{2}+2\alpha _{3}+2\alpha _{4}\}$ and

$\{\alpha _{1}+2\alpha _{2}+4\alpha _{3}+2\alpha _{4},\alpha _{1},\alpha
_{1}+2\alpha _{2}+2\alpha _{3},\alpha _{1}+2\alpha _{2}+2\alpha _{3}+2\alpha
_{4}\}$.
\end{itemize}

Hence in the maximal flag manifold the invariant subspaces are $\mathfrak{g}%
_{-\alpha }$, $\alpha $ short root, and the one dimensional subspaces
contained in $\mathfrak{g}_{-\alpha }\oplus \mathfrak{g}_{-\beta }\oplus 
\mathfrak{g}_{-\gamma }\oplus \mathfrak{g}_{-\delta }$ with $\{\alpha ,\beta
,\gamma ,\delta \}$ a $M$-equivalence class of long roots.

We will not make an extensive analysis of the other $14$ flag manifolds but
look only at the specific flag manifold $\mathbb{F}_{\Theta }$ where $\Theta
=\{\alpha _{2},\alpha _{3},\alpha _{4}\}$ which $C_{3}$ subdiagram. In this
case the $\mathfrak{z}_{\Theta }$-representation decomposes into two
irreducible components $V_{\Theta }^{1}$ and $V_{\Theta }^{2}$ with $\dim
V_{\Theta }^{1}=14$ and $\dim V_{\Theta }^{2}=1$. The sets of roots $\Pi
_{\Theta }^{i}$ corresponding to $V_{\Theta }^{i}$ are those having
coefficient $-i$ ($i=1,2$) with respect to $\alpha _{1}$. Namely, $-\Pi
_{\Theta }^{2}=\{2\alpha _{1}+3\alpha _{2}+4\alpha _{3}+2\alpha _{4}\}$ (the
highest root) and $-\Pi _{\Theta }^{1}$ contains the remaining positive
roots outside $\langle \Theta \rangle $.

Clearly the $K_{\Theta }$-representation on $V_{\Theta }^{2}$ is
irreducible. The representation on $V_{\Theta }^{1}$ is irreducible as well.
In fact, $-\alpha _{1}$ is the highest weight for the $\mathfrak{z}_{\Theta
} $-representation. Since $\langle -\alpha _{1},\alpha _{2}^{\vee }\rangle
=1 $ and $\langle -\alpha _{1},\alpha _{3}^{\vee }\rangle =\langle -\alpha
_{1},\alpha _{4}^{\vee }\rangle =0$, it follows that this is a fundamental
weight of $\mathfrak{sp}\left( 3,\mathbb{R}\right) $, namely the weight $%
\lambda _{1}+\lambda _{2}+\lambda _{3}$, where $\lambda _{i}$ has the same
meaning as in Section \ref{secirreduc}. Hence $V_{\Theta }^{1}$ is the space
of a basic representation of $\mathfrak{sp}\left( 3,\mathbb{R}\right) $. It
is known that this basic representation is irreducible by the compact
subalgebra $\mathfrak{u}\left( 3\right) $. Hence $V_{\Theta }^{1}$ is $%
K_{\Theta }$-irreducible.


\begin{thebibliography}{99}
\bibitem{br} Burstall, F. E. and J. H. Rawnsley: Twistor Theory for
Riemannian Symmetric Spaces, Springer Lect. Notes in Math.\textbf{1424}
(1990).

\bibitem{bs} Burstall, F.E. and S. Salamon: Tournaments, flags and harmonics
maps, Math. Ann.277 (1987), 249-265.

\bibitem{fh} Fulton, W. and J. Harris: Representation Theory. A first
course. Springer-Verlag.

\bibitem{he} Helgason, S.: Differential Geometry, Lie groups and Symmetric
spaces. Ac. Press (1978).

\bibitem{knp} Knapp, A.W.: Lie Groups. Beyond an Introduction. Progress in
Mathematics \textbf{140}, Birkh\"{a}user, 2004.

\bibitem{neg} Negreiros, C. J. C.: Some remarks about harmonic maps into
flag manifolds, Indiana Univ. Math. J. 37 (1988), 617-636.

\bibitem{smalg} San Martin, L.A.B.: \'{A}lgebras de Lie. Editora Unicamp
(2010).

\bibitem{smneg} San Martin, L. A. B. and C. J. C. Negreiros, Invariant
almost Hermitian structures on flag manifolds, Advances in Math., 178
(2003), 277-310.

\bibitem{smrit} San Martin, L. A. B. and R. C. J. Silva: Invariant
nearly-Kahler structures, Geom. Dedicata 121 (2006), 143-154.

\bibitem{wz} Wang, M. and W. Ziller: On normal homogeneous Einstein metrics,
Ann. Sci. Ecole norm. Sup.18(1985), 563-633.
\end{thebibliography}
\end{document}